\newcommand{\MailFile}[1]{}
\def\PaperDate{2010/01/12}
\documentclass[12pt]{article}
  \usepackage{buxmath}
  \usepackage[DinA4]{buxlayout}
  \usepackage[quotes]{bux2ref}
  \usepackage{graphicx}
  \usepackage[all]{xy}
  \usepackage{datum}
  \usepackage[german,american]{babel}
  \usepackage{float}

  \renewcommand{\qedsymbol}{\ensuremath{\Box}}

  \usepackage{email}
  \MailFile{borel.tex}
  \DoMail{hyperbolic_3_3_4.mp}
  \DoMail{w_one.mp}
  \DoMail{w_two.mp}
  \DoMail{almost_rich.mp}
  \DoMail{rich.mp}
  \DoMail{min_01.pov}
  \DoMail{min_02.pov}
  \DoMail{min_03.pov}
  \DoMail{move.mp}
  \DoNotMail{ams*, *bref, graphicx.sty, *cfg, *def, *fd}
  \DoNotMail{keyval.sty, trig.sty, xy*, babel*, *ldf}
  \DoNotMail{verbatim.sty}

  \DeclareRobustCommand{\change}[2][]{%
    \begingroup
      \def\empty{}%
      \def\arg{#2}%
      \ifx\empty\arg
        \def\next{\unskip}%
      \else
        \def\next{#2}%
      \fi
    \expandafter\endgroup\next
  }
  \let\xxx\iffalse
\xxx
  \let\oldlabel\label
  \def\label{%
    \begingroup
      \catcode`\_=11
      \dolabel
  }
  \def\dolabel#1{%
    \endgroup\oldlabel{#1}%
    \fbox{\texttt{#1}}%
  }
  \usepackage{color}
  \usepackage[normalem]{ulem}
  \usepackage{changebar}
  \DeclareRobustCommand{\delete}[1]{\begin{color}{blue}\sout{#1}\end{color}}
  \DeclareRobustCommand{\insert}[1]{\begin{color}{red}#1\end{color}}
  \DeclareRobustCommand{\swap}[2]{\delete{#1} \insert{#2}}
  \global\let\change\undefined
  \DeclareRobustCommand{\change}[2][]{%
    \cbstart
    \begingroup
      \def\next{}%
      \def\deleted{#1}%
      \def\newmat{#2}%
      \def\empty{}%
      \ifx\deleted\empty
        \def\next{\insert{#2}}%
      \else
        \ifx\newmat\empty
          \def\next{\delete{#1}}%
        \else
          \def\next{\swap{#1}{#2}}%
        \fi
      \fi
    \expandafter\endgroup\next
    \cbend
  {}}
\fi

  \InputIfFileExists{basic_notation.tex}{}{}

  \newcommand{\FType}[1][]{\ensuremath{\text{F}_{#1}}}

  \newvariable{\TheIndex}{i}
  \newvariable{\AltIndex}{j}
  \newvariable{\TheLastIndex}{k}
  \newcommand{\compose}{\circ}

  \newvariable{\Image}{\operatorname{im}}

  \newvariable{\TheGroup}{\mathcal{G}}
  \newvariable{\CovGroup}{\tilde{\TheGroup}}
  \newvariable{\TheRank}{n}
  \newvariable{\FieldOrder}{q}
  \newvariable{\TheField}{K}
  \newvariable{\AltField}{K'}
  \newvariable{\FinField}{\FFF_q}
  \newvariable{\TheIndeterminate}{t}
  \newvariable{\TheRing}{\FinFieldAd{\TheIndeterminate}}

  \newvariable{\UnivEnvel}{U}
  \newvariable{\TheMatrix}{A}
  \newvariable{\TheDatum}{\mathcal{D}}
  \newvariable{\KMalg}{L}
  \newvariable{\TheAlgebra}{\mathfrak{g}}
  \newcommand{\AlgebraicClosureOf}[1]{\overline{#1}}

  \newvariable{\GlobalField}{k}
  \newvariable{\PrimeSet}{S}
  \newvariable{\OkaSet}{\Oka_{\PrimeSet}}
  \newvariable{\AffBuilding}{X}
  \newvariable{\TwinBuilding}{\AffBuilding}
  \newvariable{\PosBuilding}{\TwinBuilding[+]}
  \newvariable{\NegBuilding}{\TwinBuilding[-]}
  \newvariable{\SubBuilding}{\TwinBuilding[+][\circ]}
  \newvariable{\HorSimplex}{\sigma_+}
  \newvariable{\HorFace}{{\HorSimplex'}}
  \newvariable{\SubSimplex}{\tilde{\tau}_+}
  \newvariable{\SubVertex}{\BarycenterOf{\HorSimplex}}
  \newvariable{\OthVertex}{\BarycenterOf{\PosSimplex}}
  \newvariable{\HorCoface}{\tau_+}
  \newvariable{\TheEnd}{\operatorname{e}}
  \newvariable{\Depth}{\operatorname{dp}}

  \global\let\Dim\undefined
  \newvariable{\Dim}{\operatorname{dim}}
  \newvariable{\TheDim}{\TheRank}
  \newvariable{\ThePlace}{p}
  \newvariable{\LocalField}{\GlobalField[\ThePlace]}
  \newvariable{\SumOfRanks}{d}
  \newvariable{\ThePrime}{p}
  \newvariable{\TheType}{m}

  \newvariable{\ArithGroup}{\TheGroupOf{\OkaSet}}
  \newvariable{\Sl}{\operatorname{SL}}

  \newvariable{\AutGroup}{G}
  \newvariable{\TheLattice}{\Gamma}
  \global\let\Stab\undefined
  \newvariable{\Stab}{\operatorname{Stab}}
  \newvariable{\TheChamber}{c}
  \newvariable{\AltChamber}{c'}
  \newvariable{\PosChamber}{\TheChamber[+]}
  \newvariable{\NegChamber}{\TheChamber[-]}
  \newvariable{\FixedChamber}{d_-}
  \newvariable{\TheApartment}{\Sigma}
  \newvariable{\TheTwinApartment}{\TheApartment}
  \newvariable{\AltTwinApartment}{\TheApartment'}
  \newvariable{\ThrTwinApartment}{\TheApartment''}
  \newvariable{\PosApartment}{\TheApartment[+]}
  \newvariable{\NegApartment}{\TheApartment[-]}
  \newvariable{\TheRetraction}{\rho}
  \newvariable{\TheIsometry}{\TheRetraction}
  \newvariable{\TheSimplex}{\tau}
  \newvariable{\AltSimplex}{\TheSimplex'}
  \newvariable{\TheEdge}{\epsilon}
  \newvariable{\TheVertex}{\nu}
  \newvariable{\AltVertex}{\nu'}
  \newvariable{\ThePoint}{x}
  \newvariable{\AltPoint}{y}
  \newvariable{\FixedPoint}{z_-}
  \newvariable{\PosPoint}{\ThePoint_+}
  \newvariable{\NegPoint}{\ThePoint_-}
  \newvariable{\op}{\operatorname{op}}

  \newvariable{\WeylGroup}{W}
  \newvariable{\TheWall}{H}
  \newvariable{\Reflection}{\sigma}
  \newcommand{\MetricCodist}[2]{\mu_*(#1,#2)}
  \newcommand{\MetricDist}[2]{\mu(#1,#2)}
  \newcommand{\CoRay}[2]{[#1,#2)}

  \newvariable{\TheHeight}{R}
  \newvariable{\TheEpsilon}{\varepsilon}
  \newvariable{\Link}{\operatorname{lk}}
  \newvariable{\Desc}{\operatorname{lk}^{\downarrow}}
  \newvariable{\SubDesc}{\operatorname{lk}^{\circ\downarrow}}
  \newvariable{\Boundary}{\partial}
  \newcommand{\ConvexHullOf}[1]{\overline{#1}}
  \newcommand{\ClosureOf}[1]{\overline{#1}}
  \newcommand{\BarycenterOf}[1]{#1[][\circ]}
  \newcommand{\hor}{\operatorname{hor}}
  \newcommand{\ver}{\operatorname{ver}}
  \newcommand{\facepart}{\partial}
  \newcommand{\cofacepart}{\delta}
  \newcommand{\join}{*}
  \newcommand{\faceof}{\leq}
  \newcommand{\strictfaceof}{<}

  \newcommand{\TheReals}{\RRR}
  \newvariable{\TheComplexes}{\CCC}
  \newcommand{\TheIntegers}{\ZZZ}
  \newcommand{\PosReals}{\RRR_{\geq 0}}
  \newvariable{\Grad}{\nabla}

  \newcommand{\isom}{\cong}
  \newvariable{\TheMorse}{h}
  \newvariable{\SubMorse}{\TheMorse^{\circ}}
  \newvariable{\Busemann}{\beta}
  \newvariable{\RoughMorse}{\TheMorse'}
  \newvariable{\DownSet}{L^{\downarrow}}
  \newvariable{\UpSet}{L^{\uparrow}}

  \newvariable{\TheResidue}{\mathcal{R}}
  \newvariable{\PosResidue}{\TheResidue[+]}
  \newvariable{\NegResidue}{\TheResidue[-]}
  \newvariable{\PosSimplex}{\TheSimplex[+]}

  \newvariable{\StdApartment}{\tilde\Sigma}
  \newvariable{\StdIsometry}{\iota}
  \newvariable{\EuclSpace}{\EEE}
  \newvariable{\DummyArg}{-}
  \newcommand{\EuclProd}[2]{\langle #1, #2 \rangle}
  \newcommand{\AngleAt}[3][]{\angle_{#1}(#2,#3)}
  \newvariable{\EuclZero}{0}
  \newvariable{\ThePolytope}{Z}
  \newvariable{\TheFace}{F}
  \newvariable{\NormalCone}{\operatorname{N}}
  \newvariable{\TheDirection}{\mathbf{n}}
  \newvariable{\TheVector}{\mathbf{z}}
  \newvariable{\AltVector}{\TheVector'}
  \newvariable{\EuclPoint}{\mathbf{x}}
  \newvariable{\EuclVertex}{\mathbf{v}}
  \newvariable{\EuclSimplex}{{\mathbf{\sigma}}}
  \newvariable{\FaceVector}{\mathbf{f}}
  \newvariable{\NormalVector}{\mathbf{n}}
  \newcommand{\NormOf}[1]{|#1|}

  \newvariable{\ZonoBase}{D}
  \newvariable{\SubBase}{D'}
  \newvariable{\ZonoVector}{\mathbf{d}}
  \newvariable{\ZonoPoint}{\mathbf{z}}
  \newvariable{\Zonotope}{\ThePolytope}
  \newvariable{\ZonoFace}{\TheFace}
  \newvariable{\RealCoeff}{\alpha}
  \newvariable{\TheZonotope}{\ThePolytope}

  \newcommand{\Away}[2][]{[#2,\infty)_{#1}}

  \global\let\Proj\undefined
  \global\let\ProjOf\undefined
  \global\let\ProjOfOf\undefined
  \global\let\ProjAd\undefined
  \global\let\ProjAdAd\undefined
  \newvariable{\Proj}{\operatorname{pr}}

  \newcommand{\goesup}{\nearrow}
  \newcommand{\goesdown}{\searrow}

  \newvariable{\TheSphere}{S}
  \newvariable{\TheSpace}{M}
  \newvariable{\NorthPole}{n}
  \newvariable{\NewPart}{A}
  \newvariable{\BadSet}{B}
  \newvariable{\FirstBadSet}{M}
  \newvariable{\SphPoint}{x}
  \newvariable{\SphComplex}{L}
  \newvariable{\SphSubcx}{K}
  \newvariable{\SphSimplex}{\TheSimplex}
  \newvariable{\DefRetr}{\rho}

  \newvariable{\SphBuilding}{\Delta}
  \newvariable{\SphChamber}{C}
  \newvariable{\SphApartment}{\Sigma}
  \newvariable{\SphConvex}{M}
  \newvariable{\SphAnticonvex}{Y}
  \newvariable{\SphRetraction}{\rho}

  \newcommand{\CAT}[1]{\ensuremath{\text{\small{CAT}}(#1)}}

\begin{document}
  \title{Finiteness Properties of Chevalley Groups over a
         Polynomial Ring over a Finite Field}
  \author{Kai-Uwe Bux \and
          Ralf Gramlich \and
          Stefan Witzel}
  \date{\datum\PaperDate}
  \maketitle
  \begin{abstract}
    \noindent
    It is known from work by H.\,Abels and P.\,Abramenko
    that for a classical $\FinField$-group $\TheGroup$ of
    rank $\TheRank$ the arithmetic lattice $\TheGroupOf{\TheRing}$
    of $\TheRing$-points is of type \FType[\TheRank-\One] provided
    that $\FieldOrder$ is large enough. We show that the statement
    is true without any assumption on $\FieldOrder$ and for any
    isotropic, absolutely
    almost simple group $\TheGroup$ defined over~$\FinField$.
  \end{abstract}

  \noindent
  Let $\GlobalField$ be a global function field and let $\TheGroup$
  be a connected, noncommutative, absolutely almost simple
  $\GlobalField$-group of positive rank. Let $\OkaSet$ be the
  ring of $\PrimeSet$-integers in $\GlobalField$. For each place
  $\ThePlace\in\PrimeSet$, there is an associated euclidean building
  $\AffBuilding[\ThePlace]$ acted upon by
  $\TheGroupOf{\LocalField}\supseteq\ArithGroup$. The
  dimension of the building $\AffBuilding[\ThePlace]$ is the
  local rank of $\TheGroup$ at the place $\ThePlace$. In
  \cite[Theorem~1.2]{Bux.Wortman:2007}, K.\,Wortman and the first
  author have shown that $\ArithGroup$ is not of type
  \FType[\SumOfRanks], where $\SumOfRanks$ is the sum of local
  ranks of $\TheGroup$ at the places in $\PrimeSet$. This settles
  the negative part of the following:
  %I don't insist on this one, it is just an older reference
  \begin{NewTh}<statement>{Rank Conjecture}[see {\cite[page 197]{Brown:1989}},
    \cite{Behr:1998}, or \cite{Bux.Wortman:2007}]
    The group $\ArithGroup$ is of type \FType[\SumOfRanks-\One]
    but not of type \FType[\SumOfRanks].
  \end{NewTh}
  Results in favor of the rank conjecture include
  \cite{Stuhler:1980} in which U.\,Stuhler shows that it holds for
  $\SlOf[\Two]{\OkaSet}$. This result has been generalized by
  Wortman and the first author in \cite{Bux.Wortman:2008} to
  arbitrary $\TheGroup$ of global rank one. Concerning higher
  ranks, H.\,Abels~\cite{Abels:1991} and
  P.\,Abramenko~\cite{Abramenko:1987} independently proved the
  rank conjecture for $\SlOf[\TheRank+\One]{\TheRing}$ provided
  that $\FieldOrder$ is sufficiently large. Abramenko has better
  bounds, but they still grow exponentially with $\TheRank$.
  In \cite{Abramenko:1996}, Abramenko has verified the rank
  conjecture for $\TheGroupOf{\TheRing}$ for classical groups
  $\TheGroup$ again under the hypothesis that $\FieldOrder$
  is sufficiently large with a bound depending only on the rank
  of $\TheGroup$. We generalize the last result.
  \begin{NewTh}<statement>{Theorem~A}
    Let $\TheGroup$ be an absolutely almost
    simple $\FinField$-group of rank $\TheRank \geq \One$.
    Then the group $\TheGroupOf{\TheRing}$
    is of type \FType[\TheRank-\One] but not of
    type \FType[\TheRank].
  \end{NewTh}
  Most of our argument will be purely geometric, and we
  shall deduce Theorem~A from:
  \begin{NewTh}<statement>{Theorem~B}
    Let $\TwinBuilding:=(\PosBuilding,\NegBuilding)$ be a
    thick, locally finite, irreducible euclidean twin building,
    and let $\AutGroup$
    be a group acting strongly
    transitively on $\TwinBuilding$
    such that any pair
    $(\PosChamber,\NegChamber)$ of chambers has a finite
    stabilizer.
    Then for an arbitrary chamber $\FixedChamber\in\NegBuilding$,
    the stabilizer
    \(
      \TheLattice
      :=
      \StabOf[\AutGroup]{\FixedChamber}
    \)
    is of type \FType[\TheDim-\One] but not of type \FType[\TheDim]
    where $\TheDim := \DimOf{\PosBuilding}$.
  \end{NewTh}
  We shall reduce Theorem~A to Theorem~B in Section~\ref{reduce}.
  We start in
  Section~\ref{heuristic} with an outline of the geometric argument
  for Theorem~B. Here, we also indicate why we have
  to confine ourselves
  to euclidean twin buildings. This is a trade off: Abramenko
  \cite{Abramenko:1996} was able to allow the compact hyperbolic
  case but had to exclude small values of $\FieldOrder$.
  The strategy proposed by H.\,Behr
  in \cite{Behr:2004} to eliminate the
  restriction on $\FieldOrder$ has not yet been carried out
  successfully.

  Only the positive parts of the claims
  are new. This is clear for Theorem~A as the fact that
  $\TheGroupOf{\TheRing}$ is not of type \FType[\TheRank] follows
  from \cite[Theorem~1.2]{Bux.Wortman:2007}. For Theorem~B,
  one can use the classification of euclidean buildings
  \cite[Chapter~28]{Weiss:2008}
  and Margulis' Arithmeticity Theorem \cite[Chapter~IX]{Margulis:1991}
  to see that $\TheLattice$ is an $\PrimeSet$-arithmetic group.
  Since $\ThePrime$-adic euclidean buildings
  do not arise in twin buildings
  \cite[Corollary~18]{Maldeghem.Muehlherr:2009},
  the group $\TheLattice$ is arithmetic over a global function field
  and \cite[Theorem~1.2]{Bux.Wortman:2007} applies.

  \paragraph{Acknowledgments.}
    We would like to thank
    Peter~Abramenko,
    Michael~Joswig,
    Bernd~Schulz, and
    Hendrik~Van~Maldeghem
    for helpful discussions and suggestions.
    Bernd~Schulz has also commented
    on an earlier version of the paper and
    suggested valuable improvements.
    We thank Bertrand~R\'emy for sharing his expertise
    on {\small RGD} systems of affine Kac-Moody groups with us.
    The first two authors also gratefully acknowledge
    the hospitality of the {\small MFO} at Oberwolfach
    allowing them to spend two weeks as
    RiP~guests
    during the final stage of the research leading to
    this article.

  \section{A geometric heuristic for Theorem~B}
  \label{heuristic}
  The method of \cite{Abramenko:1996} is to use a filtration of
  the positive partner of a twin building induced by numerical
  codistance to a fixed chamber in the negative building. We shall
  use a metric version of this idea. Our first task, therefore,
  will be to define the notion of \notion{metric codistance}
  that relates to ordinary metric distance as
  $\WeylGroup$-valued codistance relates to $\WeylGroup$-valued
  distance.

  For the following, we assume that we are given a twin
  building $\TwinBuilding=(\PosBuilding,\NegBuilding)$ that
  has a geometric realization. In particular, for any two points
  $\ThePoint[\pm]$ and $\AltPoint[\pm]$ in
  $\TwinBuilding[\pm]$, we have a metric distance between them,
  denoted by $\MetricDist{\ThePoint[\pm]}{\AltPoint[\pm]}$.
  We say for short, that $\TwinBuilding$ is a
  \notion{metric twin building}.
  The metric structure is obvious for euclidean buildings,
  where apartments are euclidean spaces. The compact hyperbolic
  case is also straightforward: here apartments have the
  geometry of hyperbolic space. In general, one could use the
  Davis realization turning the building, and each apartment,
  into a \CAT{\Zero}~space.
  Note that under any such interpretation,
  isometries in twin buildings induce isometries in the geometric
  sense.
  \begin{observation}\label{retraction_induces_isometry}
    Let $\TheTwinApartment$ and $\AltTwinApartment$ be two twin
    apartments both containing the chamber $\TheChamber$. Then,
    the retraction \cite[Exercise~5.185]{Abramenko.Brown:2008}
    \[
      \TheRetraction :=
      \TheRetraction[\TheTwinApartment,\TheChamber]
      \mapcolon \TwinBuilding \longrightarrow \TheTwinApartment
    \]
    restricts to an isometry from $\AltTwinApartment$ to
    $\TheTwinApartment$. Moreover, this isometry fixes the
    intersection $\TheTwinApartment\intersect\AltTwinApartment$
    pointwise. (Recall \cite[Exercise~5.163]{Abramenko.Brown:2008}
    that isometries of twin
    apartments also preserve codistances and in particular
    the opposition relation.)\qed
  \end{observation}
  \begin{lemma}
    Let $\TheSimplex[+]$ be a cell in $\PosBuilding$, and let
    $\TheSimplex[-]$ be a cell in $\NegBuilding$. Any two
    twin apartments $\TheTwinApartment$ and $\AltTwinApartment$
    that both contain $\TheSimplex[+]$ and $\TheSimplex[-]$ are
    isometric via an isometry fixing $\TheSimplex[+]$ and
    $\TheSimplex[-]$ pointwise.
  \end{lemma}
  \begin{proof}
    Let $\TheChamber[\pm]$ be a chamber in $\TheTwinApartment[\pm]$
    containing $\TheSimplex[\pm]$, and let $\AltChamber[\pm]$ be
    a chamber in $\AltTwinApartment[\pm]$ also containing
    $\TheSimplex[\pm]$. Let $\ThrTwinApartment$ be a twin apartment
    containing $\TheChamber[+]$ and $\AltChamber[-]$. Then
    Observation~\ref{retraction_induces_isometry} implies that
    $\ThrTwinApartment$ is isometric to
    $\TheTwinApartment$ on the one hand and to
    $\AltTwinApartment$ on the other via isometries fixing
    $\TheSimplex[+]$ and $\TheSimplex[-]$ pointwise.
  \end{proof}

  Let $\PosPoint$ be a point in $\PosBuilding$
  and $\NegPoint$ be a point in $\NegBuilding$. The
  \notion{metric codistance}
  \(
    \MetricCodist{\PosPoint}{\NegPoint}
  \)
  is the metric distance of $\PosPoint$ to the point
  $\NegPoint[][{\op[\TheTwinApartment]}]$ opposite to
  $\NegPoint$ in some twin apartment $\TheTwinApartment$
  containing $\PosPoint$ and $\NegPoint$.
  Note that the metric
  codistance does not depend on the choice of $\TheTwinApartment$,
  since any other twin apartment $\AltTwinApartment$ also
  containing $\PosPoint$ and $\NegPoint$ is isometric to
  $\TheTwinApartment$ via an isometry fixing $\PosPoint$ and
  $\NegPoint$. Since isometries respect opposition, the isometry
  takes $\NegPoint[][{\op[\TheTwinApartment]}]$ to
  $\NegPoint[][{\op[\AltTwinApartment]}]$.
  \begin{observation}\label{round_sphere}
    Let $\TheTwinApartment=(\PosApartment,\NegApartment)$ be any
    twin apartment containing $\PosPoint$ and $\NegPoint$.
    On $\PosApartment$, the metric codistance to
    $\NegPoint$ agrees with the metric distance to
    $\NegPoint[][{\op[\TheTwinApartment]}]$. In particular,
    level sets of the metric codistance to $\NegPoint$
    inside $\PosApartment$ are round spheres.\qed
  \end{observation}

  We define the \notion{geodesic ray} from $\PosPoint$ to
  $\NegPoint$ as:
  \[
    \CoRay{\PosPoint}{\NegPoint}
    :=
    \SetOf[{\AltPoint[+]\in\PosBuilding}]{
      \MetricCodist{\PosPoint}{\NegPoint}
      +
      \MetricDist{\PosPoint}{\AltPoint[+]}
      =
      \MetricCodist{\AltPoint[+]}{\NegPoint}
    }
  \]
  Rays are meaningful mostly if the metric structure
  on $\TwinBuilding$ is euclidean:
  \begin{prop}\label{rough_gradients_are_unique}
    Assume that $\TwinBuilding=(\PosBuilding,\NegBuilding)$ is
    euclidean and $\MetricCodist{\PosPoint}{\NegPoint}\neq\Zero$.
    Then, $\CoRay{\PosPoint}{\NegPoint}$ is a
    geodesic ray in the euclidean building $\PosBuilding$.
  \end{prop}
  \begin{proof}
    Assume that $\TwinBuilding$ is euclidean and that
    $\MetricCodist{\PosPoint}{\NegPoint}\neq\Zero$. Let
    $\TheTwinApartment=(\PosApartment,\NegApartment)$ be any
    twin apartment containing $\PosPoint$ and $\NegPoint$.
    Then the intersection
    \(
      \PosApartment\intersect\CoRay{\PosPoint}{\NegPoint}
    \)
    is the geodesic ray through $\PosPoint$ in $\PosApartment$
    pointing away from $\NegPoint[][{\op[\TheTwinApartment]}]$.

    Moreover, if $\AltTwinApartment$ is any other twin apartment
    containing $\PosPoint$ and $\NegPoint$, then any isometry
    from $\TheTwinApartment$ to $\AltTwinApartment$ fixing
    $\PosPoint$ and $\NegPoint$ takes the ray
    \(
      \PosApartment\intersect\CoRay{\PosPoint}{\NegPoint}
    \)
    to the ray
    \(
      \AltTwinApartment[+]\intersect\CoRay{\PosPoint}{\NegPoint}.
    \)

    Since $\CoRay{\PosPoint}{\NegPoint}$ intersects each
    twin apartment around $\PosPoint$ and $\NegPoint$ in
    a ray, the set is a union of rays issuing from $\PosPoint$.
    We have to show that there is no branching at any point
    $\PosPoint'\in\CoRay{\PosPoint}{\NegPoint}$. Note that
    if $\CoRay{\PosPoint}{\NegPoint}$ branches at $\PosPoint'$,
    then so does $\CoRay{\PosPoint'}{\NegPoint}$, so it suffices
    to show that there is no branching in the initial point.

    We want to argue by contradiction: assuming that
    branching happens, find a twin apartment containing initial
    segments of both rays. It remains to prove that
    such a twin apartment can be chosen to also contain $\NegPoint$.

    Let $\PosResidue$ be the residue in $\PosBuilding$ around
    the carrier of $\PosPoint'$,
    let $\NegChamber$ be a chamber in $\NegApartment$ containing
    $\NegPoint$, and let $\PosChamber$ be the projection of
    $\NegChamber$ into $\PosResidue$. Note that $\PosChamber$
    lies in $\PosApartment$ and contains the initial segment
    of $\CoRay{\PosPoint'}{\NegPoint}\intersect\PosApartment$.
    By \cite[Exercise~5.186]{Abramenko.Brown:2008}, for any chamber
    $\PosChamber'$ in $\PosResidue$, there is a twin apartment
    containing $\PosChamber'$, $\PosChamber$, and $\NegChamber\ni
    \NegPoint$.
    In particular, this holds for any chamber $\PosChamber'$ in
    $\PosApartment'$ containing the initial segment of
    $\CoRay{\PosPoint'}{\NegPoint}\intersect\PosApartment'$. Thus,
    we found a twin apartment, in which we could observe the
    branching of $\CoRay{\PosPoint'}{\NegPoint}$.
  \end{proof}
  \begin{observation}\label{initial_chamber}
    The proof also shows that the initial segment of
    $\CoRay{\PosPoint}{\NegPoint}$ is contained in any
    chamber of the residue $\PosResidue$ around $\PosPoint$
    that arises as a projection of a chamber containing
    $\NegPoint$.
  \end{observation}
  \begin{rem}
    Observation~\ref{rough_gradients_are_unique} holds true not
    just in euclidean twin buildings but
    in any \CAT{\Zero} twin building where apartments have the
    unique extension property for geodesic segments and all proper
    residues are spherical.
  \end{rem}

  Equipped with these geometric tools, we can make a first
  (albeit failing) attempt to prove Theorem~B. Recall the
  setup: we fix a thick, locally finite euclidean
  building $\TwinBuilding=(\PosBuilding,\NegBuilding)$,
  a group $\AutGroup$ acting strongly transitively on
  $\TwinBuilding$ such that stabilizers of twin chambers
  $(\PosChamber,\NegChamber)$ are finite, and a chamber
  $\FixedChamber$ in $\NegBuilding$. We want to determine
  the finiteness properties of the stabilizer
  \(
    \TheLattice
    :=
    \StabOf[\AutGroup]{\FixedChamber}
    .
  \)
  We also fix a point $\FixedPoint\in\FixedChamber$.

  We shall study the action of $\TheLattice$ on the euclidean
  building $\PosBuilding$. Note that the stabilizer in
  $\TheLattice$ of each chamber $\PosChamber$ in $\PosBuilding$
  is finite. Thus, all cell stabilizers
  of the $\TheLattice$-action
  are finite. As $\AutGroup$ acts strongly transitively,
  $\TheLattice$ acts transitively on the set of points
  \(
    \SetOf[\PosPoint\in\PosBuilding]{
      \MetricCodist{\PosPoint}{\FixedPoint} = \Zero
    }
    .
  \)
  For any positive real number $\TheHeight$, let
  $\PosBuildingOf{\TheHeight}$ be the maximal subcomplex of
  $\PosBuilding$ contained in the subset
  \(
    \SetOf[\PosPoint\in\PosBuilding]{
      \MetricCodist{\PosPoint}{\FixedPoint} \leq \TheHeight
    }.
  \)
  It follows from transitivity, that $\TheLattice$ acts
  cocompactly on $\PosBuildingOf{\TheHeight}$
  since $\PosBuilding$ is locally finite.

  Should it turn out that
  $\PosBuildingOf{\TheHeight}$ is $(\TheDim-\Two)$-connected for
  some $\TheHeight$,
  \cite[Proposition~1.1 and Proposition~3.1]{Brown:1987}
  would imply that $\TheLattice$ is of type
  \FType[\TheDim-\One].
  For the topological analysis, we use the Morse function
  \begin{eqnarray*}
    \RoughMorse \mapcolon \PosBuilding &
    \longrightarrow &\PosReals\\
     \PosPoint & \mapsto & \MetricCodist{\PosPoint}{\FixedPoint}
  \end{eqnarray*}
  The building $\PosBuilding$ is contractible. By standard
  arguments from combinatorial Morse theory, connectivity
  properties of sublevel complexes can be deduced from the same
  connectivity properties of descending links. It remains to argue
  that descending links are $(\TheDim-\Two)$-connected.

  Here, Proposition~\ref{rough_gradients_are_unique} is useful.
  It says that the geodesic ray $\CoRay{\PosPoint}{\FixedPoint}$
  determines a direction in $\LinkOf{\PosPoint}$, which we may
  think of as the gradient $\Grad[\PosPoint]\RoughMorse$. Because
  of Observation~\ref{round_sphere},
  the descending link at sufficiently high vertices
  should be the set of all those directions
  in $\LinkOf{\PosPoint}$ that span an obtuse angle with
  $\Grad[\PosPoint]\RoughMorse$ (\notion{gradient criterion}):
  large spheres are almost flat. Such subcomplexes of the spherical
  building $\LinkOf{\PosPoint}$ are called hemisphere complexes
  and sufficiently highly connected by results of B.\,Schulz
  \cite{Schulz:2005}.

  This strategy \emph{almost} succeeds. Generically, descending
  links are hemisphere complexes of the right dimension and
  connectivity. However, there are certain \notion{bad}
  regions in $\PosBuilding$
  (inside a single apartment, they look like corridors) where
  the descending link is not correctly detected by the gradient
  criterion. Thus, the main technical difficulty will be to perturb
  the Morse function $\RoughMorse$ so that the descending links
  inside bad corridors are improved without destroying connectivity
  of descending links in other regions.

  The remainder of this paper is organized as follows.
  After some preliminaries on zonotopes in Section~\ref{zonotopes}
  and on subcomplexes of spherical buildings in
  Section~\ref{spherical},
  we define in Section~\ref{primary} a primary Morse
  function $\TheMorse$ (the \notion{height}), which
  is $\TheLattice$-invariant and has cocompact sublevel complexes.
  The main results about this improved version
  of $\RoughMorse$ are Proposition~\ref{gradient_wellput},
  which ensures
  that gradients for $\TheMorse$ can be defined, and
  Proposition~\ref{gradient_criterion}, which says
  that the descending links with respect to $\TheMorse$ are never
  inconsistent with the gradients of $\TheMorse$. However, we cannot
  avoid that there are $\TheMorse$-horizontal edges (i.e., edges
  whose endpoints are of equal height). In order to break ties,
  we introduce a secondary and even a tertiary Morse function
  in Section~\ref{subdiv} using the \notion{depth}
  borrowed from \cite{Bux.Wortman:2008} and described in
  Section~\ref{depth}. We analyse the descending links arising
  from this Morse function in Section~\ref{down}. The final
  two sections are devoted to the proofs of Theorem~B and Theorem~A,
  respectively.

  We note that Abramenko has examples in the compact hyperbolic
  case showing that one cannot expect the analogously defined
  stabilizer $\TheLattice$ to be of type \FType[\TheRank-\One]
  in this case. Thus, it might be useful to conclude this section
  with an explanation why our strategy breaks down in the
  compact hyperbolic case.

  Roughly, we filter the building $\PosBuilding$ by metric
  codistance to $\FixedPoint$. Inside apartments, the filtration
  coincides with the filtration by metric balls centered at
  some point (opposite to $\FixedPoint$). Huge balls
  approximate horoballs. It is a feature of
  euclidean space that horospheres are hyperplanes.
  So if a huge circle runs through a vertex in a euclidean
  Coxeter complex, it will cut its star roughly in half.
  This is reason why we
  hope that, at least where the
  Morse function is large, we can expect relative links of
  our filtration to look like hemisphere complexes in spherical
  buildings.

\begin{figure}[ht]
    \null\hfill
    \includegraphics{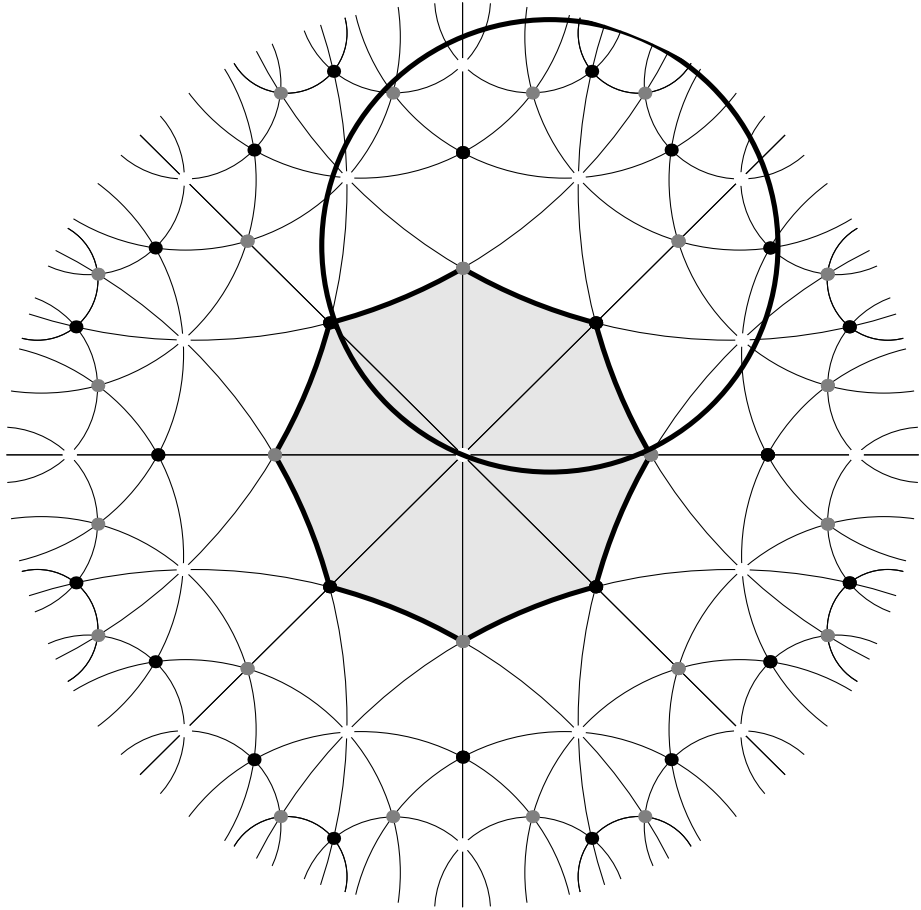}%
    \hfill\null
    \caption{a star in a hyperbolic Coxeter complex\label{fig}}
  \end{figure}

  This heuristic fails in the hyperbolic case. A horosphere
  through a vertex of a hyperbolic Coxeter complex need not split
  its star into two halves of equal size. Figure~\ref{fig} shows
  a horoball and the center vertex lies
  on its boundary circle. Only two vertices in the
  link are
  inside the horoball. This explains why relative links in
  the hyperbolic case are genuinely smaller (at least for
  filtrations based on the metric approach).

  \section{Some euclidean geometry}
  \label{zonotopes}
  Throughout this section, $\EuclSpace$ is a fixed euclidean space
  with inner product $\EuclProd{\DummyArg}{\DummyArg}$ and
  origin $\EuclZero$. Also, we fix a finite reflection group
  $\WeylGroup$ with $\EuclZero$ as a global fixed point.

  Let $\TheFace$ be a face of some (convex and compact)
  polytope $\ThePolytope\subset\EuclSpace$. The
  \notion{normal cone}
  \[
    \NormalConeOf{\TheFace}
    :=
    \SetOf[\TheDirection\in\EuclSpace]{
      \EuclProd{\TheDirection}{\TheVector}
      =
      \Max[\AltVector\in\ThePolytope]{
        \EuclProd{\TheDirection}{\AltVector}
      }
      \text{\ for all\ }
      \TheVector\in\TheFace
    }
  \]
  is the set of all $\TheDirection\in\EuclSpace$ such that the
  function $\EuclProd{\TheDirection}{\DummyArg}$ restricted
  to $\ThePolytope$ assumes its maximum on the points
  in $\TheFace$. It is a closed convex cone. We think of the
  vectors in $\NormalConeOf{\TheFace}$ as directions
  because, for any point
  $\EuclPoint\in\EuclSpace$, the closest point
  projection onto $\ThePolytope$ satisfies
  \[
    \ProjOf[\ThePolytope]{\EuclPoint}
    \in
    \TheFace
    \quad\text{\ if and only if\ }\quad
    \EuclPoint - \ProjOf[\ThePolytope]{\EuclPoint}
    \in
    \NormalConeOf{\TheFace}.
  \]
  Our reason to consider an inner product space instead of
  using a vector space and its dual is to ease phrasing of
  claims such as the following:
  \begin{lemma}\label{WeylChamber}
    Let $\ThePolytope$ be $\WeylGroup$-invariant and let $\TheFace$
    be a face of $\ThePolytope$. For any $\FaceVector\in\TheFace$
    and $\NormalVector\in\NormalConeOf{\TheFace}$, there is no wall
    with respect to $\WeylGroup$
    that separates two of the three vectors $\FaceVector$,
    $\NormalVector$, and $\FaceVector+\NormalVector$.
  \end{lemma}
  \begin{proof}
    It suffices to show that no wall
    $\TheWall$ with respect to $\WeylGroup$ separates
    $\FaceVector$ from $\NormalVector$. So assume to the
    contrary that $\TheWall$ does separate $\FaceVector$
    from $\NormalVector$. Since $\ThePolytope$
    is $\WeylGroup$-invariant, the point
    $\ReflectionOf[\TheWall]{\FaceVector}$ lies in $\ThePolytope$.
    Note that
    \(
        \ReflectionOf[\TheWall]{\FaceVector}-\FaceVector
    \)
    is orthogonal to $\TheWall$ and lies on the same side as
    $\NormalVector$ (the side opposite to $\FaceVector$). It
    follows that
    \(
      \EuclProd{\NormalVector}{
        \ReflectionOf[\TheWall]{\FaceVector}-\FaceVector
      }
      > \Zero
      ,
    \)
    whence
    \[
      \EuclProd{\NormalVector}{\ReflectionOf[\TheWall]{\FaceVector}}
      =
      \EuclProd{\NormalVector}{\FaceVector}
      +
      \EuclProd{\NormalVector}{
        \ReflectionOf[\TheWall]{\FaceVector}-\FaceVector
      }
      >
      \EuclProd{\NormalVector}{\FaceVector}
      .
   \]
   This is a contradiction as
   $\EuclProd{\NormalVector}{\FaceVector}$ is the maximum value
   of the function $\EuclProd{\NormalVector}{\DummyArg}$
   on $\ThePolytope$.
  \end{proof}

  For a finite set $\ZonoBase\subset\EuclSpace$, the convex,
  compact polytope
  \[
    \ZonotopeOf{\ZonoBase}
    :=
    \SetOf[{
      \Sum[\ZonoVector\in\ZonoBase]{
        \RealCoeff[\ZonoVector]\ZonoVector
      }
    }]{
      \Zero \leq \RealCoeff[\ZonoVector] \leq \One
      \text{\ for all\ }
      \ZonoVector\in\ZonoBase
    }
  \]
  is called the \notion{zonotope} spanned by $\ZonoBase$.
  This construction ensures:
  \begin{observation}
    Through every point $\ZonoPoint\in\ZonotopeOf{\ZonoBase}$
    and every $\ZonoVector\in\ZonoBase$, there is a line segment
    parallel to $[\EuclZero,\ZonoVector]$
    inside $\ZonotopeOf{\ZonoBase}$.\qed
  \end{observation}
  The faces of a zonotope are translated zonotopes. More
  precisely:
  \begin{lemma}
    For any direction $\NormalVector\in\EuclSpace$, the faces of
    $\ZonotopeOf{\ZonoBase}$ maximal with respect to the property
    of being orthogonal to $\NormalVector$ are translates
    of the zonotope $\ZonotopeOf{\ZonoBase[\NormalVector]}$ where
    \(
      \ZonoBase[\NormalVector]
      :=
      \SetOf[\ZonoVector\in\ZonoBase]{
        \EuclProd{\NormalVector}{\ZonoVector} = \Zero
      }
      .
    \)
    More precisely, for any point $\EuclPoint\in\EuclSpace$,
    the maximal face of $\ZonotopeOf{\ZonoBase}$ containing
    $\ProjOf[\ZonotopeOf{\ZonoBase}]{\EuclPoint}$ and
    orthogonal to
    $\NormalVector:=
    \EuclPoint-\ProjOf[\ZonotopeOf{\ZonoBase}]{\EuclPoint}$
    is given by
    \[
      \ZonoFace[\EuclPoint] :=
      \Parentheses{
      \Sum[\substack{\ZonoVector\in\ZonoBase\\
           \EuclProd{\NormalVector}{\ZonoVector} > \Zero}]{
        \ZonoVector
      }
      }
      +
      \ZonotopeOf{\ZonoBase[\NormalVector]}
      .
    \]
  \end{lemma}
  \begin{proof}
    We just observed that $\ZonoFace[\EuclPoint]$ consists precisely
    of those points in $\ZonotopeOf{\ZonoBase}$ on which the
    function $\EuclProd{\NormalVector}{\DummyArg}$ restricted
    to $\ZonotopeOf{\ZonoBase}$ is maximal. It follows that
    $\ZonoFace[\EuclPoint]$ is a face, that it is maximal among
    the faces orthogonal to $\NormalVector$, and that
    $\EuclPoint-\ProjOf[\ZonotopeOf{\ZonoBase}]{\EuclPoint}=
    \NormalVector\in\NormalConeOf{\ZonoFace[\EuclPoint]}$. The
    last statement implies
    $\ProjOf[\ZonotopeOf{\ZonoBase}]{\EuclPoint}
    \in\ZonoFace[\EuclPoint]$.
  \end{proof}

  We close this section with a consideration of the distance
  to $\TheZonotope$, i.e., the
  function $\MetricDist{\TheZonotope}{\DummyArg}$.
  As $\TheZonotope$ is convex, so is its associated distance function.
  In particular, for any
  simplex $\EuclSimplex\subset\EuclSpace$ the subset of points
  farthest away from $\TheZonotope$ contains a vertex.
  For points in $\EuclSimplex$ closest to $\TheZonotope$, we have:
  \begin{prop}\label{minimum_at_vertex}
    Let $\EuclSimplex\subset\EuclSpace$ be a simplex. Suppose
    that $\ZonoBase$ is sufficiently rich:
    $\EuclVertex-\EuclVertex'\in\ZonoBase$, for any two vertices
    $\EuclVertex$ and $\EuclVertex'$ in $\EuclSimplex$.
    Then the subset of points on $\EuclSimplex$ closest to
    $\TheZonotope:=\ZonotopeOf{\ZonoBase}$ contains a vertex.
  \end{prop}
  \begin{proof}
    Consider a point $\EuclPoint\in\EuclSimplex$ minimizing
    the distance to $\TheZonotope$.
    Without loss of generality, we may assume that
    $\EuclPoint$ lies in the relative interior of $\EuclSimplex$.
    Then,
    \(
      \NormalVector := \EuclPoint
      -
      \ProjOf[\TheZonotope]{\EuclPoint}
    \)
    is orthogonal to $\EuclSimplex$. The
    point $\ProjOf[\TheZonotope]{\EuclPoint}$
    lies in $\TheZonotope$ and by richness
    of $\ZonoBase$, there is a parallel translate
    of the zonotope $\TheZonotopeOf{\SubBase}$
    through $\ProjOf[\TheZonotope]{\EuclPoint}$
    that lies entirely in $\TheZonotope$. Here
    $\SubBase\subseteq\ZonoBase$ consists of all vectors
    representing differences of adjacent vertices in $\EuclSimplex$.
    Note
    that $\ProjOf[\TheZonotope]{\EuclPoint}$
    does not have to lie in the center of this translate.

    Nonetheless, it follows that at least one vertex
    of $\EuclSimplex$ is within distance at
    most $\NormOf{\NormalVector}$
    from this parallel copy of $\TheZonotopeOf{\SubBase}$.
    The claim follows by choice
    of $\EuclPoint$ as a point on $\EuclSimplex$ of lowest height.
  \end{proof}

  \section{Some subcomplexes of spherical buildings}
  \label{spherical}
  To deduce finiteness properties, we use the well-established
  technique of filtering a complex upon which the group acts.
  The main task, as usual, is to control the homotopy type of
  relative links that arise in the filtration.
  In this section, we collect the results concerning connectivity
  properties of those subcomplexes of spherical buildings that we
  will encounter.

  Let $\TheSpace$ be euclidean or hyperbolic space or a round
  sphere. We call an intersection of a non-empty family of closed
  half-spaces (or hemispheres in the latter case)
  \notion{demi-convex}. We call a subset
  of $\TheSpace$ \notion{fat} if it has non-empty interior.
  Note that a proper open convex subset
  of $\TheSpace$ is contained in an open hemisphere.
  \begin{observation}\label{deformation_retraction}
    Let $\NewPart\subset\TheSpace$ be fat and demi-convex
    and let $\BadSet\subset\TheSpace$ be proper, open, and convex.
    If $\NewPart$ and $\BadSet$ intersect, then
    $\NewPart\setminus\BadSet$ strongly deformation retracts
    onto the boundary part
    $\BoundaryOf{\NewPart}\setminus\BadSet$.
  \end{observation}
  \begin{proof}
    Note that $\BadSet$ intersects the interior of $\NewPart$
    since every boundary point of the convex set $\NewPart$ is
    an accumulation point of interior points because $\NewPart$
    is fat. Choose $\SphPoint$ in the intersection.
    Note that $\NewPart$ is star-like with
    regard to $\SphPoint$, and the geodesic projection away from
    $\SphPoint$ restricts to the deformation retraction we need.
  \end{proof}
  Iterated application of the same projection trick yields:
  \begin{prop}\label{push_in}
    Suppose that $\SphComplex$ is a geometric CW-complex, i.e.,
    its cells carry a spherical, euclidean, or hyperbolic
    structure in which they are demi-convex
    (i.e., each cell is an intersection of half-spaces
    in the model space).
    Let $\BadSet$ be an open subset of $\SphComplex$
    that intersects each cell in a convex set.
    Then there is a strong deformation retraction
    \[
      \DefRetr[\SphComplex]
      \mapcolon
      \SphComplex\setminus\BadSet
      \longrightarrow
%      (\SphComplex\setminus\BadSet)^\circ=:
%or has this been defined anywhere?
      \SphComplex[][\BadSet]
    \]
    of $\SphComplex\setminus\BadSet$ onto its maximal subcomplex.
  \end{prop}
  \begin{proof}
    First, we assume that $\SphComplex$ has finite dimension.
    Let $\SphSimplex$ be a maximal cell of $\SphComplex$.
    If $\TheSimplex\subseteq\BadSet$, the cell $\SphSimplex$
    does not intersect
    $\SphComplex\setminus\BadSet$ and we do not need to do anything.
    If $\SphSimplex$ avoids $\BadSet$, the map $\DefRetr$ must
    be the identity on $\SphSimplex$.
    Otherwise, let $\SphPoint$
    be a point in the intersection $\SphSimplex\intersect\BadSet$
    chosen in the relative interior of $\SphSimplex$.
    Projecting away from $\SphPoint$, as in
    Observation~\ref{deformation_retraction},
    deformation retracts $\SphSimplex\setminus\BadSet$ onto
    $\BoundaryOf{\SphSimplex}\setminus\BadSet$. The maps constructed
    for two maximal cells agree on their intersection.
    Hence we can paste all these maps together to get a deformation
    retraction of $\SphComplex\setminus\BadSet$ onto
    $\SphComplex'\setminus\BadSet$ where $\SphComplex'$ is
    $\SphComplex$ with the interiors of all maximal
    cells intersecting $\BadSet$ removed.

    Now, $\SphComplex'$ has other maximal cells, which
    might intersect $\BadSet$. Using the same construction
    for $\SphComplex'$, we obtain another deformation retraction
    $\SphComplex'\setminus\BadSet
    \rightarrow\SphComplex''\setminus\BadSet$. We keep going,
    removing more and more cells intersecting $\BadSet$.
    Since the dimension of $\SphComplex$ is finite, the process
    terminates
    after finitely many steps. The composition of the maps
    thus obtained is the strong deformation retraction from
    $\SphComplex\setminus\BadSet$ onto $\SphComplex[][\BadSet]$.
    This proves the claim for finite dimensional $\SphComplex$.

    Note that the construction is local: what it does on a cell
    is only determined by the intersection of this cell with the set
    $\BadSet$.
    Hence, the deformation retraction is compatible with
    subcomplexes. More precisely, if $\SphSubcx$ is a subcomplex of
    $\SphComplex$, then the deformation retractions
    $\DefRetr[\SphComplex]$ and $\DefRetr[\SphSubcx]$ from above
    are constructed
    such that $\DefRetr[\SphSubcx]$ is the restriction of
    $\DefRetr[\SphComplex]$ to $\SphSubcx$. It follows that the
    pair $(\SphComplex\setminus\BadSet,\SphSubcx\setminus\BadSet)$
    is homotopy equivalent
    to $(\SphComplex[][\BadSet],\SphSubcx[][\BadSet])$.
    Applying this observation to pairs of skeleta, the
    claim follows by standard arguments in the case that
    $\SphComplex$ has infinite dimension.
  \end{proof}

  Let $\SphBuilding$ be a spherical building. We regard $\SphBuilding$
  as a metric space with the angular metric. So each apartment is
  a round sphere of radius~$\One$. When $\SphBuilding$ is a
  finite building,
  the topology induced by the metric agrees with the weak topology
  it carries as a simplicial complex. For infinite buildings, both
  topologies differ and we will use the weak topology throughout
  for the building and all its subcomplexes.
  \begin{prop}\label{complement}
    Let $\SphBuilding$ be a spherical building and fix a chamber
    $\SphChamber$ in $\SphBuilding$.
    Let $\BadSet\subset\SphBuilding$
    be a subset such that, for any apartment $\SphApartment$
    containing $\SphChamber$ the intersection
    $\BadSet\intersect\SphApartment$ is a proper, open, and
    convex subset of the sphere $\SphApartment$. Then the
    space $\SphAnticonvex:=\SphBuilding\setminus\BadSet$ and its
    maximal subcomplex $\SphBuilding[][\BadSet]$ are both
    $(\DimOf{\SphBuilding}-\One)$-connected. The complex
    $\SphBuilding[][\BadSet]$ has dimension $\DimOf{\SphBuilding}$
    and hence is spherical of this dimension.
  \end{prop}
  \begin{rem}
    Using $\BadSet=\emptyset$ in Proposition~\ref{complement},
    we obtain the Solomon-Tits Theorem as a special case.
    Satz~3.5 of \cite{Schulz:2005}, whose proof
    inspired the argument given below, is the special case where
    $\BadSet$ is open, convex, and of diameter strictly less
    than $\pi$.
  \end{rem}
  \begin{proof}[of Proposition~\ref{complement}]
    We observe first that Proposition~\ref{push_in} implies that
    the subset $\SphAnticonvex$ and its maximal subcomplex
    $\SphBuilding[][\BadSet]$ are homotopy equivalent.
    Therefore, it suffices to prove that $\SphAnticonvex$ is
    $(\DimOf{\SphBuilding}-\One)$-connected.

    We have to contract spheres of dimensions up to
    $\DimOf{\SphBuilding}-\One$. Let
    $\TheSphere\subseteq\SphAnticonvex$ be such a sphere.
    Since $\TheSphere$ is compact in $\SphBuilding$, it is
    covered by a finite family of apartments and we can
    apply \cite[Lemma~3.5]{Heydebreck:2003}: there is a finite
    sequence $\SphApartment[\One],\SphApartment[\Two],\ldots,
    \SphApartment[\TheLastIndex]$ such that
    (a) each $\SphApartment[\TheIndex]$ contains $\SphChamber$,
    (b) the sphere $\TheSphere$ is contained in the union
    $\Union[\TheIndex]{\SphApartment[\TheIndex]}$, and most
    importantly,
    (c) for each $\TheIndex\geq\Two$ the intersection
    $\SphApartment[\TheIndex]\intersect
    (\SphApartment[\One]\union\cdots
    \union\SphApartment[\TheIndex-\One])$
    is a union of closed half-apartments, each of which contains
    $\SphChamber$.
    Put
    \(
      \SphComplex[\TheIndex] :=
      \SphApartment[\One] \union\cdots\union
      \SphApartment[\TheIndex]
    \)
    and observe that $\SphComplex[\TheIndex]$ is obtained
    from $\SphComplex[\TheIndex-\One]$ by gluing in the
    closure
    \(
      \NewPart[\TheIndex] :=
      \ClosureOf{
        \SphApartment[\TheIndex]\setminus
        (\SphApartment[\One]\union\cdots
        \union\SphApartment[\TheIndex-\One])
      }
    \)
    along the boundary $\BoundaryOf{\NewPart[\TheIndex]}$
    of $\NewPart[\TheIndex]$ in $\SphApartment[\TheIndex]$.
    Note that $\NewPart[\TheIndex]$ is fat and demi-convex.

    Now, we can build $\SphComplex[\TheLastIndex]\setminus\BadSet$
    inductively. We begin with $\SphComplex[\One]\setminus\BadSet$,
    which is contractible. The space
    $\SphComplex[\TheIndex]\setminus\BadSet$ is obtained
    from
    $\SphComplex[\TheIndex-\One]\setminus\BadSet$ by gluing
    in $\NewPart[\TheIndex]\setminus\BadSet$ along
    $\BoundaryOf{\NewPart[\TheIndex]}\setminus\BadSet$.
    If $\NewPart[\TheIndex]$ and $\BadSet$ are disjoint,
    this is a cellular
    extension of dimension $\DimOf{\SphBuilding}$
    as $\NewPart[\TheIndex]$ is fat. Otherwise,
    Observation~\ref{deformation_retraction} implies that
    $\NewPart[\TheIndex]\setminus\BadSet$ deformation retracts
    onto $\BoundaryOf{\NewPart[\TheIndex]}\setminus\BadSet$, whence
    $\SphComplex[\TheIndex]\setminus\BadSet$ and
    $\SphComplex[\TheIndex-\One]\setminus\BadSet$ are
    homotopy equivalent in this case. In the end,
    the sphere $\TheSphere$ can be contracted inside
    $\SphComplex[\TheLastIndex]\setminus\BadSet$.
  \end{proof}

  An interesting special case, also already noted in
  \cite{Schulz:2005}, is obtained when $\BadSet$ is chosen
  as the open $\frac{\pi}{\Two}$-ball around a fixed point
  $\NorthPole\in\SphBuilding$, which we think of as the
  \notion{north pole}. Then the complex
  $\SphBuildingOf[][\geq\frac{\pi}{\Two}]{\NorthPole} :=
  \SphBuilding[][\BadSet]$
  is a \notion{closed hemisphere complex} and
  $\DimOf{\SphBuilding}$-spherical by Proposition~\ref{complement}.
  The argument fails badly if $\BadSet$ is chosen as the
  \emph{closed ball} of radius $\frac{\pi}{\Two}$ around
  $\NorthPole$. In fact, the \notion{open hemisphere complex}
  $\SphBuildingOf[][>\frac{\pi}{\Two}]{\NorthPole}$ spanned by
  all vertices avoiding the closed ball $\BadSet$ generally is
  not $\DimOf{\SphBuilding}$-spherical: the dimension of
  $\SphBuildingOf[][>\frac{\pi}{\Two}]{\NorthPole}$ might be too
% at my oppinion
  small. The main result of Schulz is to show that this is the
  only obstruction.
  \begin{prop}[{see \cite[Page~27]{Schulz:2005}}]\label{open_hemi}
    The open hemisphere complex
    $\SphBuildingOf[][>\frac{\pi}{\Two}]{\NorthPole}$
    is spherical of dimension
    $\DimOf{\SphBuilding[\ver]}$.
    If $\SphBuilding$ is thick,
    then neither open nor closed hemisphere complexes in
    $\SphBuilding$ are contractible.
  \end{prop}
  The subcomplex $\SphBuildingOf[\ver]{\NorthPole}$ is defined as
  follows: The \notion{equator}
  $\SphBuildingOf[][=\frac{\pi}{\Two}]{\NorthPole}$ is the
  subcomplex spanned by
  those points in $\SphBuilding$ of distance $\frac{\pi}{\Two}$
  from $\NorthPole$. Recall that
  $\SphBuilding$ decomposes as a join of unique irreducible
  factors. The \notion{horizontal part}
  $\SphBuildingOf[\hor]{\NorthPole}$
  is the join of all factors fully contained in the equator.
  The complex $\SphBuildingOf[\ver]{\NorthPole}$ is the join of the
  other irreducible factors. In particular,
  \begin{equation}\label{join_dec}
    \SphBuilding = \SphBuildingOf[\hor]{\NorthPole}
    \join
    \SphBuildingOf[\ver]{\NorthPole}.
  \end{equation}

  \section{The primary Morse function}
  \label{primary}
  We now begin the proof of Theorem~B proper.
  First, let us fix a euclidean
  twin building $\TwinBuilding=(\PosBuilding,\NegBuilding)$ and
  a chamber $\FixedChamber\in\NegBuilding$. We also fix a point
  $\FixedPoint\in\FixedChamber$.

  Let $\StdApartment$ be the
  euclidean Coxeter complex upon which the apartments in
  $\PosBuilding$ are modeled. We denote by $\EuclSpace$
  the underlying euclidean space where the origin $\EuclZero$
  shall be a special vertex in $\StdApartment$. We let $\WeylGroup$
  denote the spherical Weyl group generated by the walls through
  $\EuclZero$. Finally, for this section, we choose a finite subset
  $\ZonoBase\subset\EuclSpace$. For the moment, we just require
  that it is invariant under the finite group $\WeylGroup$, but
  in the course of this work, we shall impose stronger restrictions
  upon $\ZonoBase$. The $\WeylGroup$-invariance is inherited by
  the zonotope $\TheZonotope:=\ZonotopeOf{\ZonoBase}$. In
  particular, Lemma~\ref{WeylChamber} applies.

  Consider any twin apartment
  $\TheTwinApartment=(\PosApartment,\NegApartment)$ where
  $\NegApartment$ contains the chamber $\FixedChamber$. Then,
  there is a unique point
  $\FixedPoint[][{\op[\TheTwinApartment]}]\in\PosApartment$
  opposite to $\FixedPoint$. Identifying $\PosApartment$ with
  $\EuclSpace$, we define
  \[
    \TheZonotope[\TheTwinApartment]
    :=
    \FixedPoint[][{\op[\TheTwinApartment]}]
    +
    \TheZonotope
  \]
  Observe that $\TheZonotope[\TheTwinApartment]$ is well-defined
  since any identification of $\PosApartment$ with $\EuclSpace$
  that respects the structure of the underlying labeled Coxeter
  complexes gives rise to the same
  set $\TheZonotope[\TheTwinApartment]$ since $\TheZonotope$
  is $\WeylGroup$-invariant. The same consideration shows:
  \begin{observation}\label{zonotopes_match}
    Let $\AltTwinApartment=(\PosApartment',\NegApartment')$ be
    another twin apartment containing $\FixedChamber$. Then
    the isometry
    \(
      \TheIsometry
      =
      \TheRetraction[\AltTwinApartment,\FixedChamber]
      \mapcolon \TheTwinApartment \rightarrow
      \AltTwinApartment
    \)
    from Observation~\ref{retraction_induces_isometry} takes
    $\TheZonotope[\TheTwinApartment]$ to
    $\TheZonotope[\AltTwinApartment]$.\qed
  \end{observation}
  \begin{rem}
    One can do even better. The positive partners in twin
    apartments as above are parameterized by the chambers in
    $\PosBuilding$ opposite to $\FixedChamber$. It follows that
    the isometries fixing $\FixedChamber$ from
    Observation~\ref{retraction_induces_isometry}
    are canonical and allow one to identify all such
    apartments with $\StdApartment$ in a compatible way, i.e.,
    there is a familily of isometries
    \(
      \StdIsometry[\TheTwinApartment]\mapcolon\StdApartment
      \rightarrow\PosApartment
    \)
    such that the diagrams
    \[
      \xymatrix{
        {\PosApartment}
          \ar[rr]^{\TheRetraction[\AltTwinApartment,\FixedChamber]}
        &
        &
        {\PosApartment'}
        \\
        &
        {\StdApartment}
          \ar[ul]^{{\StdIsometry[\TheTwinApartment]}}
          \ar[ur]_{{\StdIsometry[\AltTwinApartment]}}
        &
        \\
      }
    \]
    all commute. Finally, identifying $\StdApartment$ and
    $\EuclSpace$, the zonotope $\TheZonotope[\TheTwinApartment]$
    would be
    well-defined even for zonotopes $\TheZonotope\subset\EuclSpace$
    that are not $\WeylGroup$-invariant.
  \end{rem}

  Given any point $\PosPoint\in\PosBuilding$, we chose a twin
  apartment $\TheTwinApartment=(\PosBuilding,\NegBuilding)$
  containing $\PosPoint$ and $\FixedChamber$. We define
  the \notion{height} of $\PosPoint$ to be the metric distance
  \[
    \TheMorseOf{\PosPoint}
    :=
    \MetricDist{\TheZonotope[\TheTwinApartment]}{\PosPoint}
  \]
  from $\PosPoint$ to the convex, compact polytope
  $\TheZonotope[\TheTwinApartment]$.
  Observation~\ref{zonotopes_match} implies that
  $\TheMorseOf{\PosPoint}$ is independent of the chosen
  twin apartment $\TheTwinApartment$.
  \begin{observation}\label{convex}
    Let $\TheTwinApartment=(\PosApartment,\NegApartment)$ be
    a twin apartment containing $\FixedChamber$.
    The restriction of $\TheMorse$ to $\PosApartment$ is a convex
    function as the metric distance to a convex compact
    polytope. In particular, at least one highest point on
    any simplex is a vertex.\qed
  \end{observation}

  Turning to gradients of $\TheMorse$,
  we can first define the ray
  $\Away[\TheTwinApartment]{\PosPoint}$
  relative to the
  twin apartment $\TheTwinApartment$ as the direction of the
  geodesic ray in $\PosApartment$ through $\PosPoint$ away
  from $\ProjOf[{\TheZonotope[\TheTwinApartment]}]{\PosPoint}$.
  Here, we assume that $\TheMorseOf{\PosPoint}>\Zero$.
  \begin{prop}\label{gradient_wellput}
    Let $\PosPoint\in\PosApartment$ be a point with
    $\TheMorseOf{\PosPoint}>\Zero$.
    Let $\TheTwinApartment=(\PosApartment,\NegApartment)$
    and $\AltTwinApartment=(\PosApartment',\NegApartment')$
    be two twin apartments both containing $\PosPoint$ and
    $\FixedChamber$. Then the two geodesic rays
    \(
      \Away[\TheTwinApartment]{\PosPoint}
    \)
    and
    \(
      \Away[\AltTwinApartment]{\PosPoint}
    \)
    coincide.
  \end{prop}
  \begin{proof}
    Let $\PosChamber$ be the projection of $\FixedChamber$
    into the residue around the carrier of $\PosPoint$. Since
    $\TheTwinApartment$ contains
    $\PosPoint$ and $\FixedChamber$, it also contains
    $\PosChamber$. By Observation~\ref{initial_chamber}, the
    gradient
% In Observation it stated very precisely, that the gradient lies "in *every*
% chamber s.th."; as soon as there are two, of course it cannot be carried
% by any of them
    $\Grad[\PosPoint]{\RoughMorse}$ is contained in $\PosChamber$.

    We identify $\PosApartment$ with $\EuclSpace$ so
    that $\TheZonotope[\TheTwinApartment]\subset\PosApartment$
    corresponds to $\TheZonotope\subset\EuclSpace$. In particular,
    $\op[\TheTwinApartment]{\FixedPoint}$ corresponds
    to $\EuclZero$. Let $\EuclPoint\in\EuclSpace$ be the vector
    corresponding to $\PosPoint\in\PosApartment$. Then
    $\FaceVector:=\ProjOf[\TheZonotope]{\EuclPoint}$ corresponds
    to $\ProjOf[{\TheZonotope[\TheTwinApartment]}]{\PosPoint}$.
    Put $\NormalVector:=\EuclPoint-\FaceVector$. By
    Lemma~\ref{WeylChamber}, any $\WeylGroup$-chamber containing
    $\EuclPoint=\NormalVector+\FaceVector$ also contains
    $\NormalVector$. %See Figure~\ref{grad_comp}.
    Note that the vector $\EuclPoint$ is parallel to
    $\Grad[\PosPoint]{\RoughMorse}$. Hence, the chamber
    $\PosChamber$ also contains an initial segment of
    $\Away[\TheTwinApartment]{\PosPoint}$, which is
    parallel to $\NormalVector$.
    \begin{figure}[ht]
      \null\hfill
      \includegraphics{w_one.ps}\hfill
      \includegraphics{w_two.ps}\hfill\null
      \caption{comparing $\Grad[\PosPoint]{\RoughMorse}$ and
               $\Grad[\PosPoint]{\TheMorse}$\label{grad_comp}}
      \small In each figure, the circle is a level set
      for $\RoughMorse$, and the perpendicular arrow indicates
      $\Grad[\PosPoint]{\RoughMorse}$.
      The other closed curve is a level set for $\TheMorse$. The
      $\WeylGroup$-chamber for the dihedral group of order $8$
      based at $\PosPoint$ and
      containing $\Grad[\PosPoint]{\RoughMorse}$ also contains
      $\Grad[\PosPoint]{\TheMorse}$. Consequently, any chamber
      of the underlying Coxeter complex supporting
      $\Grad[\PosPoint]{\RoughMorse}$ also carries
      $\Grad[\PosPoint]{\TheMorse}$. Note that $\PosPoint$ does
      not need to be a vertex.
    \end{figure}

    The same argument applies to $\AltTwinApartment$. Hence
    both apartments contain $\PosChamber$ and at least the initial
    segments of $\Away[\TheTwinApartment]{\PosPoint}$ and
    $\Away[\AltTwinApartment]{\PosPoint}$ agree.

    The intersection of the two rays is a closed set (in each
    of the rays) and for
    any point $\AltPoint[+]$ in the intersection, we have:
    \[
      \Away[\TheTwinApartment]{\AltPoint[+]}
      \subseteq
      \Away[\TheTwinApartment]{\PosPoint}
      \qquad\text{and}\qquad
      \Away[\AltTwinApartment]{\AltPoint[+]}
      \subseteq
      \Away[\AltTwinApartment]{\PosPoint}
    \]
    Therefore, the argument from above also shows that the rays
    $\Away[\TheTwinApartment]{\PosPoint}$ and
    $\Away[\AltTwinApartment]{\PosPoint}$ share a little open
    segment around $\AltPoint[+]$. Thus, the intersection of
    the two rays is open in each of the rays. Since each ray
    is connected the claim follows.
  \end{proof}
  Proposition~\ref{gradient_wellput} implies that we can define the
  \notion{flow line}
  \(
    \Away{\PosPoint}
  \)
  through $\PosPoint$ as the geodesic ray
  \(
    \Away[\TheTwinApartment]{\PosPoint}
  \)
  in any twin apartment containing $\PosPoint$.
  We define the
  \notion{gradient}
  \(
    \Grad[\PosPoint]\TheMorse
  \)
  as the direction of
  $\Away{\PosPoint}$ at $\PosPoint$.

  A good deal of our analysis regards the interplay of the
  simplicial structure on $\PosBuilding$ and the height $\TheMorse$.
  We start with the following:
  \begin{observation}\label{orthogonal_gradient}
    Let $\PosSimplex$ be a simplex in $\PosBuilding$ and let
    $\PosPoint$ be a point in $\PosSimplex$. If
    $\Grad[\PosPoint]\TheMorse$ is orthogonal to $\PosSimplex$,
    then $\PosPoint$ is a point where $\TheMorse$ assumes its
    minimum value on $\PosSimplex$.
  \end{observation}
  \begin{proof}
    We choose a twin apartment
    $\TheTwinApartment=(\PosApartment,\NegApartment)$ containing
    $\PosSimplex$ and $\FixedChamber$. In $\PosApartment$,
    $\TheMorse$ is given as the function
    $\MetricDist{\TheZonotope[\TheTwinApartment]}{\DummyArg}$.
    The claim now follows since $\TheZonotope[\TheTwinApartment]$
    is convex and $\PosSimplex$ spans an affine subspace.
  \end{proof}

  Recall that $\EuclSpace$ is identified with the model apartment
  $\StdApartment$.
  We say that a finite subset $\ZonoBase\subset\EuclSpace$
  is \notion{almost rich} if it is $\WeylGroup$-invariant and
  contains all differences of adjacent vertices in $\StdApartment$.
  If $\ZonoBase$ is almost rich, $\TheMorse$ does not suffer from
  the obvious deficiencies of $\RoughMorse$.
  \begin{figure}[ht]
    \null\hfill
    \includegraphics{almost_rich.ps}\hfill\null
    \caption{an almost rich zonotope
    for $\tilde{\text{\textsf{B}}}_2$\label{almost_rich_zono}}
    \small
    This figure illustrates some $\TheMorse$-level sets. The lightly
    shaded areas and the white corridors are the normal cones
    for the zonotope.
  \end{figure}
  \begin{prop}\label{gradient_criterion}
    Assume that $\ZonoBase$ is almost rich.
    Let $\TheEdge$ be an edge in $\PosBuilding$ connecting
    the vertices $\TheVertex$ and $\AltVertex$. Then the following
    hold:
    \begin{enumerate}
      \item
        The function $\TheMorse$ is monotonic on $\TheEdge$.
      \item
        The angle
        $\AngleAt[\TheVertex]{\TheEdge}{\Grad[\TheVertex]\TheMorse}
        >\frac{\pi}{\Two}$ if and only if
        $\TheMorseOf{\TheVertex}>\TheMorseOf{\AltVertex}$.
    \end{enumerate}
  \end{prop}
  \begin{proof}
    The function $\TheMorse$ is convex and attains its
    minimum at a boundary point by
    Proposition~\ref{minimum_at_vertex}.
    This proves the first
    claim. If
    \(
      \AngleAt[\TheVertex]{\TheEdge}{\Grad[\TheVertex]\TheMorse}
      \neq\frac{\pi}{\Two}
    \)
    the height $\TheMorse$ changes when one moves from the vertex
    $\TheVertex$ infinitesimally into the edge $\TheEdge$. If the
    angle is obtuse, the height decreases; if the angle is
    acute, the height increases. Accordingly, $\TheVertex$ must
    be the highest or lowest point, respectively, as $\TheMorse$
    is monotonic on $\TheEdge$.
    Observation~\ref{orthogonal_gradient} covers the remaining
    case that
    $\Grad[\TheVertex]\TheMorse$ is orthogonal to $\TheEdge$.
  \end{proof}

  \section{Simplices of constant height}
  \label{h_horizontal}
  A simplex $\PosSimplex$ in $\PosBuilding$ is
  \notion{$\TheMorse$-horizontal} if $\TheMorse$ restricts
  to a constant function on $\PosSimplex$.
  \begin{observation}\label{flow_on_horizontal}
    Let $\PosSimplex$ be an $\TheMorse$-horizontal simplex.
    Then all flow lines issuing in $\PosSimplex$ are
    pairwise parallel and orthogonal to $\PosSimplex$.
  \end{observation}
  \begin{proof}
    The claim is clear
    for points in the relative interior of $\PosSimplex$.
    It follows for points on the boundary by continuity.
  \end{proof}
  It follows that we can talk about the \notion{gradient}
  $\Grad[\PosSimplex]\TheMorse$ for an $\TheMorse$-horizontal
  simplices $\PosSimplex$, which we identify with the gradient
  at the barycenter of $\PosSimplex$.

  Let $\PosPoint$ be a point carried by the simplex $\PosSimplex$
  in $\PosBuilding$. We think of the
  \notion{link} $\LinkOf{\PosPoint}$ as the space of directions
  issuing from $\PosPoint$. The link of $\LinkOf{\PosSimplex}$
  is the space of directions at $\PosPoint$ orthogonal to
  $\PosSimplex$. It does not depend on the particular point
  $\PosPoint$ carried by $\PosSimplex$. Both links are
  spherical buildings and can be regarded as
  metric spaces via the angular metric. The point link splits
  as a spherical join
  \begin{equation}\label{point_link}
    \LinkOf{\PosPoint}
    =
    \BoundaryOf{\PosSimplex}
    \join
    \LinkOf{\PosSimplex}
  \end{equation}
  where $\BoundaryOf{\PosSimplex}$ is the round sphere of directions
  at $\PosPoint$ that do not leave $\PosSimplex$. Note that
  $\BoundaryOf{\PosSimplex}$ has an obvious simplicial structure
  being the boundary of a simplex. The spherical building
  $\LinkOf{\PosSimplex}$ also has a simplicial structure,
  whose face poset corresponds to the poset of cofaces
  of $\PosSimplex$.

  Now we specialize to the case where $\HorSimplex$ is a
  horizontal simplex in $\PosBuilding$.
  For concreteness, we think of $\LinkOf{\HorSimplex}$ as centered
  at the barycenter
  $\BarycenterOf{\HorSimplex}$ of $\HorSimplex$.
  Note that $\Grad[\BarycenterOf{\HorSimplex}]\TheMorse
  \in\LinkOf{\HorSimplex}$ as
  a consequence of Observation~\ref{flow_on_horizontal}. We regard
  the distinguished point
  $\Grad[\BarycenterOf{\HorSimplex}]\TheMorse$ as the
  north pole in the spherical building $\LinkOf{\HorSimplex}$.
  Thus for any horizontal simplex, the link decomposes as
  in~(\ref{join_dec}):
  \begin{equation}\label{simplex_link}
    \LinkOf{\HorSimplex}
    =
    \LinkOf[\hor]{\HorSimplex}
    \join
    \LinkOf[\ver]{\HorSimplex}
  \end{equation}
  of the link into the horizontal and vertical part of
  $\LinkOf{\HorSimplex}$ with respect to the north pole
  $\Grad[\BarycenterOf{\HorSimplex}]\TheMorse$. We call
  the horizontal part
  $\LinkOf[\hor]{\HorSimplex}$ the \notion{horizontal link},
  and we call the vertical part
  $\LinkOf[\ver]{\HorSimplex}$ the \notion{vertical link}
  of $\HorSimplex$.
  Beware that the vertical link can contain equatorial
  simplices; and consequently, not every $\TheMorse$-horizontal
  coface of $\HorSimplex$ defines a simplex
  in $\LinkOf[\hor]{\HorSimplex}$.

  \section{The depth of horizontal simplices}
  \label{depth}
  Horizontal simplices are the main obstacle for the analysis of
  the cocompact filtration of $\PosBuilding$ by height. We will
  use the method of \cite{Bux.Wortman:2008} to cope with this
  difficulty. Here, we mostly
  follow \cite[Section~5]{Bux.Wortman:2008}.

  Let $\HorSimplex$ be an $\TheMorse$-horizontal simplex in
  $\PosBuilding$. By Observation~\ref{flow_on_horizontal}, the
  flow lines starting in $\HorSimplex$ are pairwise parallel
  geodesic
  rays in $\PosBuilding$ and therefore, they define a point
  $\TheEndOf{\HorSimplex}$
  in the spherical building at infinity. Let $\Busemann$ be
  a Busemann function centered at that point. Since the
  flow lines are orthogonal to $\HorSimplex$, the function
  $\Busemann$ is constant on $\HorSimplex$, i.e., the simplex
  $\HorSimplex$ is \notion{$\Busemann$-horizontal}.
  The notion of the
  horizontal and vertical link of $\HorSimplex$ defined above
  agree with the notions in \cite[Section~5]{Bux.Wortman:2008}, whence
  we can use some results therein directly.
  \begin{lemma}\label{equivalence}
    For any $\TheMorse$-horizontal simplex $\HorSimplex$, there
    is a unique face $\HorSimplex[][\min]$ such that for any proper
    face $\HorFace<\HorSimplex$, the following equivalence holds
    \[
      \HorSimplex\setminus\HorFace
      \in
      \LinkOf[\hor]{\HorFace}
      \quad\text{\ if and only if\ }\quad
      \HorSimplex[][\min]
      \faceof
      \HorFace
      .
    \]
  \end{lemma}
  \begin{proof}
    The statement is \cite[Lemma~5.2]{Bux.Wortman:2008} for
    $\Busemann$-horizontal simplices. Since we can find a Busemann
    function $\Busemann$, for which $\HorSimplex$ and all its
    faces are $\Busemann$-horizontal, the claim follows.
  \end{proof}
  In the same way, the following lemma is an immediate consequence
  of \cite[Observation~5.3]{Bux.Wortman:2008}.
  \begin{lemma}\label{wedging}
    Suppose $\HorSimplex[][\min]\faceof\HorFace\faceof\HorSimplex$.
    Then $\HorSimplex[][\min]=\HorFace[][\min]$.\qed
  \end{lemma}
  \begin{figure}[ht]
    \null\hfill
    \includegraphics[scale=1]{min_01.ps}\hfill
    \includegraphics[scale=1]{min_02.ps}\hfill\null
    \caption{the face $\HorSimplex[][\min]$\label{min_illustrated}}
    \small
    Both figures take place inside the Coxeter
    complex~$\tilde{\text{\textsf{B}}}_3$. In the picture on
    the left hand side,
    the black vertex is the face $\HorSimplex[][\min]$ of the
    horizontal solidly colored $\Two$-simplex $\HorSimplex$.
    The two edges of $\HorSimplex$ containing $\HorSimplex[][\min]$
    illustrate
    Lemma~\ref{wedging}. In the picture on the right, the
    horizontal simplex $\HorSimplex$ is the center edge. Here, we
    have $\HorSimplex=\HorSimplex[][\min]$.
  \end{figure}

  For any two $\TheMorse$-horizontal simplices $\HorSimplex$
  and $\HorFace$, we define \notion{going up} as
  \[
    \HorFace \goesup \HorSimplex
    \quad
    :\Longleftrightarrow
    \quad
    \HorFace=\HorSimplex[][\min] \neq \HorSimplex
  \]
  and \notion{going down} as
  \[
    \HorSimplex \goesdown \HorFace
    \quad
    :\Longleftrightarrow
    \quad
    \HorSimplex[][\min] \not\faceof \HorFace
    \strictfaceof \HorSimplex.
  \]
  We define a \notion{move} as either going up or going down.
  \begin{figure}[ht]
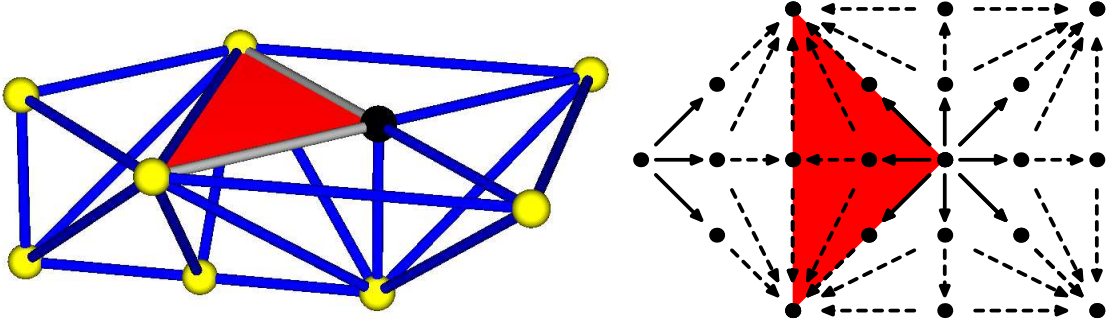

    \null\hfill
    \includegraphics{min_03.ps}\hfill
    \includegraphics{move.ps}\hfill\null
    \caption{sequences of moves\label{move_sequence}}
    \small
    This figure continues Figure~\ref{min_illustrated}; it also
    takes place inside a Coxeter complex of
    type~$\tilde{\text{\textsf{B}}}_3$. It shows a possible
    patch of horizontal $\Two$-simplices. Each dot in the
    picture on the right represents
    a simplex; for orientation, one horizontal $\Two$-simplex
    has been filled in. Arrows indicate moves: solid arrows
    represent going up, whereas dashed arrows represent going
    down. Note that there are no moves between triangles and
    short edges.
  \end{figure}
  \begin{observation}\label{move}
    If there is a move from $\HorSimplex$ to $\HorSimplex'$,
    then either $\HorSimplex$ is a face of $\HorSimplex'$ or
    $\HorSimplex'$ is a face of $\HorSimplex$. In
    either case, we have
    $\TheEndOf{\HorSimplex}=\TheEndOf{\HorSimplex'}$.\qed
  \end{observation}
  \begin{prop}
    There is a uniform upper bound, depending only on the
    building $\PosBuilding$, on the length of any sequence
    of moves.
  \end{prop}
  \begin{proof}
    Consider a sequence of moves starting at the horizontal
    simplex $\HorSimplex$. By Observation~\ref{move}, for any
    simplex $\HorSimplex'$ encountered along that sequence,
    we have $\TheEndOf{\HorSimplex}=\TheEndOf{\HorSimplex'}$.
    Let $\Busemann$ be a Busemann function centered at
    $\TheEndOf{\HorSimplex}$. It follows that the given
    sequence of moves is a ``$\Busemann$-sequence'' consisting
    of ``$\Busemann$-moves'' as considered
    in \cite[Proposition~5.4]{Bux.Wortman:2008}, where
    the existence of a uniform bound (depending only on the
    dimension of $\PosBuilding$) on the length of any such
    sequence is proved.
  \end{proof}

  We define the \notion{depth} $\DepthOf{\HorSimplex}$
  of
  an $\TheMorse$-horizontal simplex $\HorSimplex$
  as the maximum length of a sequence of moves
  starting at $\HorSimplex$.
  \begin{rem}
    Since not every $\Busemann$-move is a legal $\TheMorse$-move,
    the depth as defined here will generally be lower than the
    depth used in \cite{Bux.Wortman:2008}.
  \end{rem}

  \section{Subdividing along horizontal simplices}
  \label{subdiv}
  In this section, we assume that the
  zonotope $\TheZonotope$ is defined
  via an almost rich set $\ZonoBase$. Then, being connected by
  an $\TheMorse$-horizontal edge is an equivalence relation on
  the vertices of a given simplex $\PosSimplex$ and the equivalence
  classes correspond to the  maximal horizontal faces
  of $\PosSimplex$. This fact allows us to mimic the subdivision
  rule used in \cite[Section~6]{Bux.Wortman:2008}.

  Let $\SubBuilding$ be the simplicial subdivision of $\PosBuilding$
  whose vertices are precisely the
  barycenters $\BarycenterOf{\HorSimplex}$ of
  $\TheMorse$-horizontal simplices $\HorSimplex$.
  More precisely, we subdivide
  each horizontal simplex barycentrically; any simplex is the
  simplicial join of its maximal horizontal faces and carries the
  induced subdivision. Note that this rule of subdividing is
  compatible with inclusion of faces.
  \begin{observation}\label{flag_complex}
    The building $\PosBuilding$ is a flag complex, and so is
    the subdivision $\SubBuilding$.\qed
  \end{observation}

  As in \cite[Observation~6.1]{Bux.Wortman:2008}, the link
  $\LinkOf{\SubVertex}$ of a vertex $\SubVertex\in\SubBuilding$
  corresponding to the horizontal simplex $\HorSimplex$
  decomposes as a join
  \begin{equation}\label{coarse_decomposition}
    \LinkOf{\SubVertex}
    =
    \LinkOf[\facepart]{\SubVertex}
    \join
    \LinkOf[\cofacepart]{\SubVertex}
  \end{equation}
  where $\LinkOf[\facepart]{\SubVertex}$ is the barycentric
  subdivision of $\BoundaryOf{\HorSimplex}$ and
  $\LinkOf[\cofacepart]{\SubVertex}$ is the induced subdivision of
  $\LinkOf{\HorSimplex}$. The latter again decomposes as the join of
  its horizontal and vertical parts (see~\ref{simplex_link});
  however, this decomposition is
  not (immediately) compatible with the simplicial structure on
  $\LinkOf[\cofacepart]{\SubVertex}$.

  The building $\PosBuilding$ carries the geometric structure
  of a euclidean simplicial complex. In particular, barycenters
  have a geometric meaning and simplices in $\SubBuilding$ could
  be regarded as honest subsets of simplices in $\PosBuilding$.
  Regarded this way, $\TheMorse$ already \emph{is} a function on
  $\SubBuilding$. However,
  we only use $\TheMorse$ to define values on the vertices in
  $\SubBuilding$, which is not a problem at all since vertices
  in $\SubBuilding$ correspond to simplices in $\PosBuilding$
  on which $\TheMorse$ is already constant.

  We use the following Morse function on the vertices of
  $\SubBuilding$:
  \begin{eqnarray*}
    \SubMorse \mapcolon {\SubBuilding}^{(\Zero)}
    & \longrightarrow &
    \TheReals\times\TheReals\times\TheReals
    \\
    \SubVertex & \mapsto &
    \TupelOf{
      \TheMorseOf{\HorSimplex},
      \DepthOf{\HorSimplex},
      \DimOf{\HorSimplex}
    }
  \end{eqnarray*}
  We use lexicographic
  comparison to order $\TheReals\times\TheReals\times\TheReals$.
  \begin{observation}\label{no_critical_edges}
    There are no $\SubMorse$-horizontal edges, i.e.:
    if $\SubVertex$ and $\OthVertex$ are adjacent in
    $\SubBuilding$, then
    $\SubMorseOf{\SubVertex}\neq\SubMorseOf{\OthVertex}$.
  \end{observation}
  \begin{proof}
    If there is an edge between $\SubVertex$ and $\OthVertex$ and
    $\TheMorseOf{\HorSimplex}=\TheMorseOf{\PosSimplex}$, then
    $\HorSimplex$ is a face of $\PosSimplex$ or vice versa. In
    either case, the dimensions differ.
  \end{proof}

  \section{Descending links}
  \label{down}
  Let $\SubVertex$ be a vertex in $\SubBuilding$. The
  \notion{descending link} $\DescOf{\SubVertex}$ is the
  subcomplex of $\LinkOf{\SubVertex}$ spanned by all neighbors
  of $\SubVertex$ of strictly lower $\SubMorse$-height.
  In this section, we shall determine the
  connectivity of descending links. The argument follows
  closely the blueprint in \cite[Section~6]{Bux.Wortman:2008}.
  \begin{observation}\label{desc_coarse_join}
    Links in flag complexes are flag complexes. Hence the
    decomposition~(\ref{coarse_decomposition}) induces a
    decomposition
    \begin{equation}\label{coarse_descending}
      \DescOf{\SubVertex}
      =
      \DescOf[\facepart]{\SubVertex}
      \join
      \DescOf[\cofacepart]{\SubVertex}
    \end{equation}
    of the descending link. Here
    \(
      \DescOf[\facepart]{\SubVertex}
      :=
      \LinkOf[\facepart]{\SubVertex} \intersect \DescOf{\SubVertex}
     \)
    and
    \(
      \DescOf[\cofacepart]{\SubVertex}
      :=
      \LinkOf[\cofacepart]{\SubVertex}
      \intersect \DescOf{\SubVertex}.
    \)\qed
  \end{observation}
  \begin{lemma}\label{i_am_not_my_min}
    For any horizontal simplex $\HorSimplex$ with
    $\HorSimplex\neq\HorSimplex[][\min]$, the descending link
    $\DescOf{\SubVertex}$ is contractible.
  \end{lemma}
  \begin{proof}
    The argument is the same as
    in \cite[Lemma~6.5]{Bux.Wortman:2008}. We reproduce the main
    steps for the convenience of the reader.

    Because of the decomposition~(\ref{coarse_descending}), it
    suffices to show that $\DescOf[\facepart]{\SubVertex}$ is
    contractible. Recall that $\LinkOf[\facepart]{\SubVertex}$
    is the barycentric subdivision of $\BoundaryOf{\HorSimplex}$.
    As $\HorSimplex$ is horizontal, $\TheMorse$ will not decide
    among its faces whether they define descending directions.

    As $\HorSimplex[][\min]\neq\HorSimplex$, one can go up
    $\HorSimplex[][\min]\goesup\HorSimplex$, whence
    $\DepthOf{\HorSimplex[][\min]}>\DepthOf{\HorSimplex}$.
    It follows that $\HorSimplex[][\min]$ does not define a
    vertex in $\DescOf[\facepart]{\SubVertex}$.

    For any face $\HorFace\strictfaceof\HorSimplex$ with
    $\HorSimplex[][\min]\not\faceof\HorFace$, one can go down
    $\HorSimplex\goesdown\HorFace$. Hence
    $\DepthOf{\HorSimplex}>\DepthOf{\HorFace}$ and the barycenter
    of $\HorFace$ belongs to $\DescOf{\SubVertex}$.

    We do not know what happens to barycenters of faces
    $\HorFace$ with
    $\HorSimplex[][\min]\strictfaceof\HorFace$. Nonetheless,
    the face part $\LinkOf[\facepart]{\SubVertex}$ is a sphere.
    Its descending part is obtained by
    puncturing the sphere at the barycenter
    of $\HorSimplex[][\min]$. The
    cofaces of $\HorSimplex[][\min]$ may or may not be
    non-descending, but that only determines the size of the
    puncture: $\HorSimplex[][\min]$ will provide a cone point
    for the hole; and its complement, the descending face
    part $\DescOf[\facepart]{\SubVertex}$, is contractible.
  \end{proof}

  Now we turn to the descending links of vertices
  $\SubVertex$ with $\HorSimplex[][\min]=\HorSimplex$.
  We begin with the face part.
  \begin{lemma}\label{face_part_descending}
    For any horizontal simplex $\HorSimplex$ with
    $\HorSimplex=\HorSimplex[][\min]$, the face part
    $\LinkOf[\facepart]{\SubVertex}$ is completely descending.
  \end{lemma}
  \begin{proof}
    For any proper face $\HorFace\strictfaceof\HorSimplex$,
    one goes down $\HorSimplex\goesdown\HorFace$ since
    $\HorSimplex[][\min]=\HorSimplex\not\faceof
    \HorFace\strictfaceof\HorSimplex$. Thus,
    $\DepthOf{\HorSimplex}>\DepthOf{\HorFace}$.
  \end{proof}

  The coface part is more difficult. Ignoring some subdivision
  issues for the moment, it decomposes as the join of its vertical
  and horizontal part. It turns out that the depth of simplices
  behaves oppositely in both regions.
  \begin{lemma}\label{vertical_depth}
    Let $\HorSimplex$ be an $\TheMorse$-horizontal simplex and
    let $\HorCoface$ be an $\TheMorse$-horizontal
    coface of $\HorSimplex$,
    i.e., assume $\HorSimplex\strictfaceof\HorCoface$.
    If $\HorSimplex[][\min]=\HorSimplex$
    and $\HorCoface\setminus\HorSimplex$ is a simplex not
    completely contained in
    the horizontal part $\LinkOf[\hor]{\HorSimplex}$, then
    $\DepthOf{\HorCoface}>\DepthOf{\HorSimplex}$. In particular,
    the conclusion holds if $\HorCoface\setminus\HorSimplex$
    lies in the vertical link $\LinkOf[\ver]{\HorSimplex}$.
  \end{lemma}
  \begin{proof}
    By Lemma~\ref{equivalence}, we have
    $\HorCoface[][\min]\not\faceof
    \HorSimplex\strictfaceof\HorCoface$. Hence we go down
    $\HorCoface\goesdown\HorSimplex$, whence
    $\DepthOf{\HorCoface}>\DepthOf{\HorSimplex}$.
  \end{proof}
  \begin{lemma}\label{horizontal_depth}
    Let $\HorSimplex$ be an $\TheMorse$-horizontal simplex and
    let $\HorCoface$ be an $\TheMorse$-horizontal
    coface of $\HorSimplex$,
    i.e., assume $\HorSimplex\strictfaceof\HorCoface$.
    If $\HorSimplex[][\min]=\HorSimplex$
    and $\HorCoface\setminus\HorSimplex$ is a simplex
    completely contained in
    the horizontal part $\LinkOf[\hor]{\HorSimplex}$, then
    $\DepthOf{\HorSimplex}>\DepthOf{\HorCoface}$.
  \end{lemma}
  \begin{proof}
    By Lemma~\ref{equivalence}, we have
    $\HorCoface[][\min]\faceof
    \HorSimplex\strictfaceof\HorCoface$. By Lemma~\ref{wedging},
    we conclude $\HorCoface[][\min]=\HorSimplex[][\min]$.
    Hence we go up
    $\HorSimplex\goesup\HorCoface$, whence
    $\DepthOf{\HorSimplex}>\DepthOf{\HorCoface}$.
  \end{proof}
  To summarize: in the vertical link, the depth is always biased
  toward being ascending; in the horizontal link, the depth is
  always biased in favor of descent.

  This also helps with the subdivision issues alluded to above.
  The link $\LinkOf[\cofacepart]{\SubVertex}$ is a subdivision
  (inherited from $\SubBuilding$)
  of $\LinkOf{\HorSimplex}$, and the latter decomposes as a join
  of its horizontal and vertical part. This decomposition is
  \emph{not} compatible with the subdivision. The problem is that
  there can be an $\TheMorse$-horizontal coface $\HorCoface$ of
  $\HorSimplex$ that has vertices in the horizontal and vertical
  part. The barycentric subdivision of $\HorCoface$ does not
  respect the join decomposition. However, Lemma~\ref{vertical_depth}
  implies that for such $\HorCoface$, the barycenter
  $\BarycenterOf{\HorCoface}$ cannot belong to
  $\DescOf[\cofacepart]{\SubVertex}$. Thus, we see the
  following:
  \begin{lemma}\label{down_join}
    The coface part $\DescOf[\cofacepart]{\SubVertex}$
    of the descending link decomposes as a join:
    \[
      \DescOf[\cofacepart]{\SubVertex}
      =
      \SubDescOf[\hor]{\HorSimplex}
      \join
      \SubDescOf[\ver]{\HorSimplex}
    \]
    Here $\SubDescOf[\hor]{\HorSimplex}$ is the subdivision of
    $\DescOf[\hor]{\HorSimplex}$ obtained from barycentrically
    subdividing $\TheMorse$-horizontal simplices;
    and $\SubDescOf[\ver]{\HorSimplex}$ is defined mutatis mutandis.\qed
  \end{lemma}

  Knowing the depth component of $\SubMorse$, we turn to the
  height component. It behaves nicely on the vertical component.
  \begin{lemma}\label{coface_vertical}
    Let $\HorSimplex$ be an $\TheMorse$-horizontal simplex. The
    \notion{down set}
    \[
      \DownSet
      :=
      \SetOf[{
        \TheVertex\in\LinkOf{\HorSimplex}
      }]{
        \TheMorseOf{\TheVertex}<\TheMorseOf{\HorSimplex}
      }
    \]
    of strictly lower vertices in its link spans the open
    hemisphere complex in $\LinkOf[\ver]{\HorSimplex}$ with
    respect to the north
    pole $\Grad[\BarycenterOf{\HorSimplex}]\TheMorse$.
  \end{lemma}
  \begin{proof}
    First note that a vertex $\TheVertex\in\LinkOf{\HorSimplex}$
    below $\HorSimplex$ lies in the
    vertical link $\LinkOf{\HorSimplex}$: if it was
    horizontal, Observation~\ref{orthogonal_gradient} would
    rule out $\TheMorseOf{\TheVertex}<\TheMorseOf{\HorSimplex}$.

    Now let $\TheVertex$ be any vertex
    of $\LinkOf[\ver]{\HorSimplex}$. Fix a vertex $\AltVertex$
    in $\HorSimplex$ and let $\TheEdge$ be the edge from
    $\AltVertex$ to $\TheVertex$.
    Using Proposition~\ref{gradient_criterion} and the standing
    assumption that $\ZonoBase$ is almost rich, we have
    \[
      \AngleAt[\AltVertex]{\TheEdge}{\Grad[\AltVertex]\TheMorse}
      >\frac{\pi}{\Two}
      \quad
      \text{\ if and only if\ }
      \quad
      \TheVertex\in\DownSet
      .
    \]
    Now note that $\Grad[\AltVertex]\TheMorse$ and
    $\Grad[\BarycenterOf{\HorSimplex}]\TheMorse$ are
    parallel and perpendicular to the line connecting
    $\AltVertex$ to the barycenter $\BarycenterOf{\HorSimplex}$.
    Thus, if $\TheEdge'$ is the straight line segment
    connecting $\BarycenterOf{\HorSimplex}$ to $\TheVertex$,
    we have
    \[
      \AngleAt[\BarycenterOf{\HorSimplex}]{\TheEdge'
        }{\Grad[\BarycenterOf{\HorSimplex}]\TheMorse}
      >\frac{\pi}{\Two}
      \quad
      \text{\ if and only if\ }
      \quad
      \TheVertex\in\DownSet
      .
    \]
    Hence $\DownSet$ is the vertex set of the open
    hemisphere complex with respect to the north pole
    $\Grad[\BarycenterOf{\HorSimplex}]\TheMorse$ in
    $\LinkOf{\HorSimplex}$.
  \end{proof}

  It remains to study the height $\TheMorse$ on horizontal
  links $\LinkOf[\hor]{\HorSimplex}$. Recall that the
  $\WeylGroup$-invariant subset $\ZonoBase\subset\EuclSpace$
  is almost rich, i.e., it contains all vectors connecting
  adjacent vertices in $\StdApartment\isom\EuclSpace$.
  We call $\ZonoBase$ \notion{rich} if it contains
  each vector connecting two vertices whose closed stars intersect.
  \begin{figure}[ht]
    \null\hfill
    \includegraphics{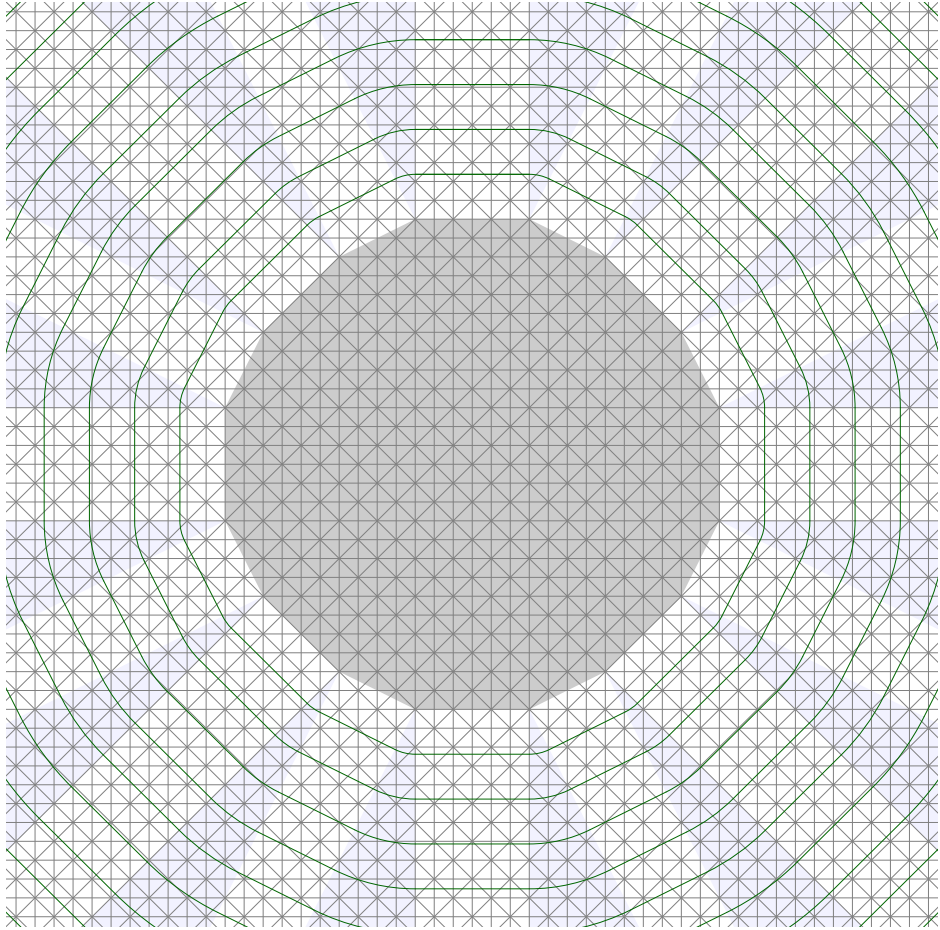}\hfill\null
    \caption{a rich zonotope
    for $\tilde{\text{\textsf{B}}}_2$\label{rich_zonotope}}
    \small This figure illustrates some $\TheMorse$-level sets.
    The lightly shaded areas
    and the white corridors are the normal cones for the zonotope.
  \end{figure}
  \begin{lemma}\label{rich}
    Assume that $\ZonoBase$ is rich. Let
    $\TheTwinApartment=(\PosApartment,\NegApartment)$ be a
    twin apartment containing $\FixedChamber$, let
    $\PosSimplex$ be a simplex in $\PosApartment$, and let
    $\TheHeight$ be the maximum value of $\TheMorse$
    on $\PosSimplex$. Then, there is a vertex $\TheVertex$
    in the set
    \(
      \UpSet :=
      \SetOf[\AltVertex\in\LinkOf{\PosSimplex}]{
        \TheMorseOf{\AltVertex} > \TheHeight
      }
    \)
    of all vertices in the link $\LinkOf{\PosSimplex}$ strictly
    higher than $\PosSimplex$ such that $\TheMorseOf{\TheVertex}$
    is the minimum value of $\TheMorse$ on the convex hull
    of $\UpSet$. In particular, the convex hull
    $\ConvexHullOf{\UpSet}$
    is disjoint
    from $\PosSimplex$ and therefore $\ConvexHullOf{\UpSet}$ and
    $\PosSimplex$ are separated by a hyperplane in $\PosApartment$.
  \end{lemma}
  \begin{proof}
    By Proposition~\ref{minimum_at_vertex}, any simplex spanned
    by vertices in $\UpSet$ has a vertex where $\TheMorse$
    assumes its minimum value on that simplex. Since these
    simplices cover the convex hull of $\UpSet$ (Carath\'eordory's
    Theorem), the claim follows.
  \end{proof}
  \begin{rem}
    One could also include in $\ZonoBase$ difference vectors
    arising from barycenters of simplices with intersecting stars.
    This would slightly simplify the proof of Lemma~\ref{coface_vertical}.
    (We suggest the name \notion{filthy rich} for such $\ZonoBase$.)
  \end{rem}

  Finally, we can paste the pieces together.
  \begin{lemma}\label{i_am_my_min}
    Assume that $\ZonoBase$ is rich.
    Let $\HorSimplex$ an $\TheMorse$-horizontal simplex with
    $\HorSimplex[][\min]=\HorSimplex$. Then the descending link
    $\DescOf{\SubVertex}$ of the barycenter
    $\SubVertex$ is spherical of dimension $\DimOf{\PosBuilding}-\One$.
  \end{lemma}
  \begin{proof}
    By Observation~\ref{desc_coarse_join} and
    Lemma~\ref{down_join}:
    \[
      \DescOf{\SubVertex}
      =
      \DescOf[\facepart]{\SubVertex}
      \join
      \SubDescOf[\hor]{\HorSimplex}
      \join
      \SubDescOf[\ver]{\HorSimplex}
    \]
    The face part $\DescOf[\facepart]{\SubVertex}$ is a sphere
    by Lemma~\ref{face_part_descending}.
    The vertical coface part
    \(
      \SubDescOf[\ver]{\HorSimplex}
    \)
    is a subdivided open hemisphere complex in $\LinkOf{\HorSimplex}$
    by Lemma~\ref{coface_vertical} and Lemma~\ref{vertical_depth}.
    By Proposition~\ref{open_hemi}, this part is spherical of dimension
    $\DimOf{\LinkOf[\ver]{\HorSimplex}}$.

    It remains to consider the horizontal
    part $\SubDescOf[\hor]{\HorSimplex}$.
    Ultimately, we want
    to apply Proposition~\ref{complement}. So first, we regard
    $\LinkOf{\HorSimplex}$ as a residue in $\PosBuilding$.
    Let $\PosChamber$ be the projection of $\FixedChamber$ in
    the residue $\LinkOf{\HorSimplex}$.
    Let $\TheTwinApartment=(\PosApartment,\NegApartment)$ be
    a twin apartment containing $\PosChamber$ and $\FixedChamber$.
    Note that $\PosChamber$ contains
    $\Grad[\HorSimplex]\RoughMorse$ and
    $\Grad[\HorSimplex]\TheMorse$.

    Let $\UpSet$ be the set of all vertices in
    $\LinkOf{\HorSimplex}\intersect\PosApartment$ strictly
    higher than $\HorSimplex$. Lemma~\ref{rich} implies that
    the subcomplex of the spherical apartment
    $\LinkOf{\HorSimplex}\intersect\PosApartment$
    spanned by $\UpSet$
    is convex and has diameter strictly
    less than $\pi$. It follows that there is an open, proper,
    convex subset $\FirstBadSet$ of the round sphere
    $\LinkOf{\HorSimplex}\intersect\PosApartment$ that contains
    all vertices in $\UpSet$ and no other
    vertices of the apartment. E.g., one could
    take $\FirstBadSet$ to be the convex hull
    of sufficiently small balls around the points in
    $\UpSet$.

    Any other twin apartment $\AltTwinApartment$ that contains
    $\PosChamber$ and $\FixedChamber$ is isometric to
    $\TheTwinApartment$ via an isometry leaving $\PosChamber$
    fixed. Define $\FirstBadSet[\AltTwinApartment]$ to be
    the isometric image
    of $\FirstBadSet$ in
    $\LinkOf{\HorSimplex}\intersect\AltTwinApartment[+]$.
    As $\TheMorse$ is compatible with the isometry,
    $\FirstBadSet[\AltTwinApartment]$ is an open, proper, convex subset containing
    precisely those vertices of
    $\LinkOf{\HorSimplex}\intersect\AltTwinApartment[+]$
    that are strictly higher than $\HorSimplex$.

    Now put
    \[
      \BadSet :=
      \LinkOf[\hor]{\HorSimplex}
      \intersect
      \Union[\TheTwinApartment]{
        \FirstBadSet[\TheTwinApartment]
      }
    \]
    where $\TheTwinApartment$ ranges over all twin apartments
    containing $\PosChamber$ and $\FixedChamber$.
    Since $\PosChamber$ is the projection of $\FixedChamber$ into
    the residue $\LinkOf{\HorSimplex}$, any apartment of
    $\LinkOf[\hor]{\HorSimplex}$ that contains the chamber
    $\PosChamber\intersect\LinkOf[\hor]{\HorSimplex}$
    comes from such a twin apartment. Thus,
    $\BadSet$ is an open subset of $\LinkOf[\hor]{\HorSimplex}$
    that satisfies the hypotheses of Proposition~\ref{complement}.
    It follows that the subcomplex of $\LinkOf[\hor]{\HorSimplex}$
    spanned by all those vertices below or at the height
    of $\HorSimplex$ is spherical.
    By Lemma~\ref{horizontal_depth}, the descending horizontal
    link $\SubDescOf[\hor]{\HorSimplex}$ is a subdivision of
    precisely this spherical complex.
  \end{proof}
  Combining Lemma~\ref{i_am_not_my_min} and Lemma~\ref{i_am_my_min},
  we see:
  \begin{prop}\label{descent_spherical}
    Assume that $\ZonoBase$ is rich. The descending link
    $\DescOf{\SubVertex}$
    of any vertex $\SubVertex\in\SubBuilding$ is spherical of
    dimension $\DimOf{\PosBuilding}-\One$.\qed
  \end{prop}

  Of course, the generic vertex will not have neighbors of equal
  height. It is only along some regions that we encounter strange
  links. In the generic case, the descending link is a hemisphere
  complex (in the generic case, open and closed
  makes no difference). Using Proposition~\ref{open_hemi} for
  thick buildings, we conclude:
  \begin{observation}\label{non_contractible}
    There are arbitrarily high vertices with non-contractible
    descending links.\qed
  \end{observation}

  \section{Finiteness properties: proof of Theorem~B}
  \label{proof}
  Finally, we assume that the twin building $\TwinBuilding$ is
  locally finite. We also assume that the set $\ZonoBase$,
  defining the zonotope, is $\WeylGroup$-invariant and rich.
  E.g., one could chose $\ZonoBase$ to consist precisely of the
  difference vectors of any pair of vertices
  in $\StdApartment\isom\EuclSpace$ whose closed stars intersect.
  \begin{observation}
    The group $\TheLattice$ acts simplicially on $\SubBuilding$.
    Simplex stabilizers of the action are finite.\qed
  \end{observation}
  \begin{observation}
    The function $\SubMorse$ is $\TheLattice$-invariant, and
    its sublevel complexes are $\TheLattice$-cocompact.\qed
  \end{observation}
  \begin{proof}[of Theorem~B]
    Given the topological properties of descending links, the
    deduction of finiteness properties is routine.

    Since $\TheLattice$ acts cocompactly, there are only finitely
    many $\TheLattice$-orbits of vertices in $\SubBuilding$
    below any given
    $\SubMorse$-bound in $\TheReals\times\TheReals\times\TheReals$.
    In particular, only finitely many elements in
    $\TheReals\times\TheReals\times\TheReals$ arise as
    values of $\SubMorse$ below any given bound.
    Define
    \(
      \SubBuildingOf{\TheIndex}
    \)
    to be the subcomplex of $\SubBuilding$ spanned by all vertices
    $\SubVertex$
    such that there are at most $\TheIndex$ values in the image
    $\ImageOf{\SubMorse}$
    that are strictly below $\SubMorseOf{\SubVertex}$.

    By Observation~\ref{no_critical_edges},
    there are no $\SubMorse$-horizontal edges in
    $\SubBuilding$. Thus,
    $\SubBuildingOf{\TheIndex+\One}\setminus
    \SubBuildingOf{\TheIndex}$ does not contain adjacent
    vertices.

    Recall that $\TheDim$ denotes the
    dimension $\DimOf{\PosBuilding}$.
    For any vertex $\SubVertex\in
    \SubBuildingOf{\TheIndex+\One}\setminus
    \SubBuildingOf{\TheIndex}$ the relative link
    $\LinkOf{\SubVertex}\intersect\SubBuildingOf{\TheIndex}$
    is precisely the descending link $\DescOf{\SubVertex}$.
    By Proposition~\ref{descent_spherical},
    descending links are spherical of dimension
    $\DimOf{\PosBuilding}-\One$. Thus, the complex
    $\SubBuildingOf{\TheIndex+\One}$ is obtained up to
    homotopy equivalence from $\SubBuildingOf{\TheIndex}$ by
    attaching $\TheDim$-cells.
    Observation~\ref{non_contractible} ensures that at
    infinitely many stages the extension is non-trivial.

    The group $\TheLattice$ acts on $\SubBuilding$ with
    finite stabilizers by hypothesis. Thus, all hypotheses
    of the criterion
    \cite[Corollary~3.3]{Brown:1987} are
    satisfied and $\TheLattice$
    is of type \FType[\TheDim-\One] but not of type
    \FType[\TheDim].
  \end{proof}

  \section{Deducing Theorem~A from Theorem~B}
  \label{reduce}
  The gap between Theorem~A and Theorem~B is bridged by the
  construction of a twin building for the group
  \(
      \TheGroupOf{\FinFieldAd{
        \TheIndeterminate,\TheIndeterminate[][-\One]
      }}
      .
  \)
  Although certainly known to the experts, we were not able to
  find a clean reference. For this reason, we outline the
  construction for connected, simply connected groups (this will be
  sufficient for the application to finiteness properties).
% c.f. P. Abramenko, H. Abels and the like; note also that a mentioning
% in the acknowledgements does not count as first mentioning, as can
% be seen by the example of B. Schulz
%  B.\,R\'emy is going to provide a citable reference in the near
%  future.

  First, we deal with split groups.
  \begin{prop}\label{functor}
    Let $\TheField$ be a field of arbitrary characteristic and
    let $\TheGroup$ be an isotropic, connected,
    simply connected, almost simple,
    split $\TheField$-group. Then the functor
    \(
      \TheGroupOf{\DummyArgAd{
        \TheIndeterminate,\TheIndeterminate[][-\One]
      }}
    \)
    is a Kac-Moody functor.
  \end{prop}
  We should explain our notation:
  to any field $\AltField$, the functor above assigns
  the group of
  \(
    \AltFieldAd{\TheIndeterminate,\TheIndeterminate[][-\One]}
  \)-points of $\TheGroup$.
  \begin{proof}
    By \cite[Theorem~16.3.2]{Springer:1998} and
    \cite[\S II]{Chevalley:1955}, an isotropic, connected,
    simply connected, almost simple
    $\TheField$-group that splits over $\TheField$ is a
    Chevalley group. It follows that the group scheme
    $\TheGroupOf{\DummyArg}$ is defined over $\TheIntegers$.
    Hence the functor
    \(
      \TheGroupOf{\DummyArgAd{
        \TheIndeterminate,\TheIndeterminate[][-\One]
      }}
    \)
    can be defined for all fields.

    A Kac-Moody functor is associated to a root
    datum $\TheDatum$, the main part
    of which is a generalized Cartan matrix $\TheMatrix$.
    Classically, this kind of datum classifies reductive groups
    over the complex numbers. There, the generalized
    Cartan matrix
    is not really generalized and defines a
    finite Coxeter group. Kac-Moody functors were
    defined by Tits \cite{Tits:1987} in the case where the
    generalized Cartan matrix defines an arbitrary Coxeter group.

    In order to recognize
    \(
      \TheGroupOf{\DummyArgAd{
        \TheIndeterminate,\TheIndeterminate[][-\One]
      }}
    \)
    as a Kac-Moody functor, we have to correctly identify its
    defining datum $\TheDatum$. Since the group $\TheGroup$ is
    simply connected, we only have to choose the generalized
    Cartan
    matrix $\TheMatrix$. Here, we use the unique generalized
    Cartan matrix
    given by a euclidean Coxeter diagram extending the spherical
    diagram as defined by $\TheGroup$.

    To show that
    \(
      \TheGroupOf{\DummyArgAd{
        \TheIndeterminate,\TheIndeterminate[][-\One]
      }}
    \)
    is the Kac-Moody functor associated to $\TheDatum$, one needs
    to verify the axioms {\small (KMG~1)} through
    {\small (KMG~9)} in \cite{Tits:1987}. All axioms are straight
    forward to check; however {\small (KMG~5)} and
    {\small (KMG~6)} involve the complex Kac-Moody algebra
    $\KMalgOf{\TheMatrix}$ associated to the given Cartan matrix.
    To verify these, one needs to know that
    $\KMalgOf{\TheMatrix}$ is the universal central extension
    of the Lie algebra
    \(
      \TheAlgebraOf{\TheComplexesAd{
        \TheIndeterminate,\TheIndeterminate[][-\One]
      }}
    \)
    where $\TheAlgebra$ is the Lie algebra associated
    to $\TheGroup$. See e.g., \cite[Theorem~9.11]{Kac:1990} or
    \cite[Section~5.2]{Pressley.Segal:1986}.
  \end{proof}

% c.f. P. Abramenko, H. Abels and the like; note also that a mentioning
% in the acknowledgements does not count as first mentioning, as can
% be seen by the example of B. Schulz
  In \cite{Remy:2002}, B.\,R\'emy has extended
  the construction
  to non-split groups using the method of Galois descent.
  \begin{prop}\label{non_split}
    Let $\TheGroup$ be an isotropic, connected,
    simply connected, almost simple
    group defined over the finite field $\FinField$. Then
    the functor
    \(
      \TheGroupOf{\DummyArgAd{
        \TheIndeterminate,\TheIndeterminate[][-\One]
      }}
    \)
    is an almost split $\FinField$-form of a Kac-Moody group
    defined over the algebraic
    closure $\AlgebraicClosureOf{\FinField}$.
  \end{prop}
  \begin{proof*}
    First, $\TheGroup$ splits
    over $\AlgebraicClosureOf{\FinField}$. Hence,
    \(
      \TheGroupOf{\DummyArgAd{
        \TheIndeterminate,\TheIndeterminate[][-\One]
      }}
    \)
    is a Kac-Moody functor over $\AlgebraicClosureOf{\FinField}$
    by the preceding proposition. Let $\TheDatum$ be the
    associated root datum.

    Note that the conditions {\small(KMG 6)} through
    {\small(KMG 9)} ensure that the ``abstract'' and ``constructive''
    Kac-Moody functors associated to $\TheDatum$ coincide
    \cite[Theorem 1']{Tits:1987},
    which holds in particular for
    \(
      \TheGroupOf{\DummyArgAd{
        \TheIndeterminate,\TheIndeterminate[][-\One]
      }}
      .
    \)
    This is relevant as R\'emy discusses Galois descent
    for constructive Kac-Moody functors.

    The claim follows from \cite[Section~11]{Remy:2002}
    once a list of conditions scattered throughout that section
    have been verified. Checking individual axioms is
    straightforward,
    the hard part (left to the reader) is making sure that
    no condition is overlooked. Here is the list:
    \begin{description}
      \item[\normalfont{\small(PREALG 1)} {[p.~257]}]
        One needs to know that $\UnivEnvel[\TheDatum]$ is the
        $\TheIntegers$-form of the universal enveloping algebra
        of $\KMalgOf{\TheMatrix}$. Its $\FinField$-form
        is obtained by the Galois action.
      \item[\normalfont{\small(PREALG 2)} {[p.~257]}]
        Straightforward.
      \item[\normalfont{\small(SGR)} {[p.~266]}]
        Straightforward.
      \item[\normalfont{\small(ALG 1)} {[p.~267]}]
        Use Definition~11.2.1 on page~261.
      \item[\normalfont{\small(ALG 2)} {[p.~267]}]
        Straightforward.
      \item[\normalfont{\small(PRD)} {[p.~273]}]
        Observe that the Galois group acts trivially on
        $\TheIndeterminate$ and $\TheIndeterminate[][-\One]$.\qed
    \end{description}
  \end{proof*}

  We are finally closing in on twin buildings.
  \begin{prop}
    Let $\TheGroup$ be as in Proposition~\ref{non_split}.
    The group
    \(
      \TheGroupOf{\FinFieldAd{
        \TheIndeterminate,\TheIndeterminate[][-\One]
      }}
    \)
    has an {\small RGD}~system.
  \end{prop}
  \begin{proof*}
    This follows from \cite[Theorem~12.4.3]{Remy:2002}; but
    once again, we need to verify hypotheses. This time, we have
    to deal with only two:
    \begin{description}
      \item[\normalfont{{\small (DCS$_1$)} [p.~284]}]
        This holds as $\TheGroup$ splits already over a finite
        field extension of $\FinField$.
      \item[\normalfont{{\small (DCS$_2$)} [p.~284]}]
        This follows from $\FinField$ being a finite, and
        hence perfect field.\qed
    \end{description}
  \end{proof*}
  \begin{prop}\label{twin_is_there}
    Let $\TheGroup$ be an isotropic, connected,
    simply connected, almost
    simple group defined over the finite field $\FinField$
    (i.e., $\TheGroup$ is as in Proposition~\ref{non_split}).
    Then there is a thick, locally finite, irreducible
    euclidean twin building
    $\TwinBuilding=(\PosBuilding,\NegBuilding)$ on which
    \(
      \TheGroupOf{\FinFieldAd{
        \TheIndeterminate,\TheIndeterminate[][-\One]
      }}
    \)
    acts strongly transitively. Moreover, $\PosBuilding$
    and $\NegBuilding$ are
    \(
      \TheGroupOf{\FinFieldAd{
        \TheIndeterminate,\TheIndeterminate[][-\One]
      }}
    \)-equivariantly
    isomorphic to the euclidean
    building associated to
    $\TheGroupOf{\FinFieldOfOf{\TheIndeterminate}}$
    and
    $\TheGroupOf{\FinFieldOfOf{\TheIndeterminate[][-\One]}}$,
    respectively.
  \end{prop}
  \begin{proof}
    By the preceding proposition, the group
    \(
      \TheGroupOf{\FinFieldAd{
        \TheIndeterminate,\TheIndeterminate[][-\One]
      }}
    \)
    has an {\small RGD}~system. By
    \cite[Theorem~8.80 and Theorem~8.81]{Abramenko.Brown:2008},
    we find an associated twin building upon which the group
    acts strongly transitively. Theorem~8.81 also tells us
    that the root groups act simply transitively, which implies
    that the twin building is thick and locally finite. That it is
    irreducible and euclidean is clear as we chose the generalized
    Cartan matrix $\TheMatrix$ back in
    in the proof Proposition~\ref{functor} to match the spherical
    type of $\TheGroup$, which is almost simple.

    The identification of $\PosBuilding$ and $\NegBuilding$
    with the euclidean building associated to
    $\TheGroupOf{\FinFieldOfOf{\TheIndeterminate}}$,
    respectively
    $\TheGroupOf{\FinFieldOfOf{\TheIndeterminate[][-\One]}}$,
    follows from the functoriality of the construction
    \cite[5.1.2]{Rousseau:1977}.
  \end{proof}
  \begin{rem}
    It also follows from \cite[Theorem~8.81]{Abramenko.Brown:2008}
    that the building thus constructed is Moufang.
  \end{rem}

  \begin{rem}
    For split groups, Abramenko gives the {\small RGD}~system
    explicitly in \cite[Example~3, page~18]{Abramenko:1996}.
    He also derives {\small RGD}~systems for groups of the types
    ${}^2\widetilde{\text{\textsf{A}}}_{\TheRank}$ and
    ${}^2\widetilde{\text{\textsf{D}}}_{\TheRank}$
    in \cite[Chapter~III.1]{Abramenko:1996}. Hence, the only
    types not covered by his explicit computations are
    ${}^3\widetilde{\text{\textsf{D}}}_4$ and
    ${}^2\widetilde{\text{\textsf{E}}}_6$. The marginal gain
    also explains why we merely sketched the general argument.
  \end{rem}

  \begin{proof}[of Theorem~A using Theorem~B]
    Let $\TheGroup$ be as in Theorem~A, i.e., $\TheGroup$
    is an isotropic, almost simple group defined over
    the finite field $\FinField$. We may assume that $\TheGroup$
    is connected since the connected component of the identity
    element has finite index.

    Let $\CovGroup$ be its
    ``universal cover'', i.e., a simply connected, isotropic,
    almost simple $\FinField$-group which allows for a central
    isogeny onto $\TheGroup$. By \cite[Satz~2]{Behr:1968}, the
    image of $\CovGroupOf{\FinFieldAd{\TheIndeterminate}}$
    in $\TheGroupOf{\FinFieldAd{\TheIndeterminate}}$ under the
    isogeny has finite index. As the isogeny has finite kernel,
    the finiteness properties of
    \(
      \TheGroupOf{\FinFieldAd{\TheIndeterminate}}
    \)
    and
    \(
      \CovGroupOf{\FinFieldAd{\TheIndeterminate}}
    \)
    coincide. Hence we may assume without loss of generality
    that $\TheGroup$ is simply connected.

    Now, we can apply Proposition~\ref{twin_is_there}. Hence
    there is a thick, locally finite, irreducible euclidean
    twin building
    $\TwinBuilding=(\PosBuilding,\NegBuilding)$ on which
    \(
      \AutGroup:=\TheGroupOf{
        \FinFieldAd{
          \TheIndeterminate,\TheIndeterminate[][-\One]
        }
      }
    \)
    acts strongly transitively. The $\AutGroup$-equivariant
    isomorphisms of $\PosBuilding$ and $\NegBuilding$ to the
    euclidean building associated to
    \(
      \TheGroupOf{
        \FinFieldOfOf{\TheIndeterminate}
      }
    \)
    implies first that the $\FinField$-rank $\TheRank$ of
    $\TheGroup$ is the dimension of the building.
    It also implies
    that stabilizers of pairs $(\PosChamber,\NegChamber)$
    of chambers are finite: they are compact and discrete.
% Is the following statement obvious?
    Finally, the group
    $\TheGroupOf{\FinFieldAd{\TheIndeterminate}}$ is
    commensurable to the stabilizer
    $\TheLattice := \StabOf[\AutGroup]{\FixedChamber}$ for some
    chamber $\FixedChamber\in\NegBuilding$. Since finiteness
    properties are invariant under commensurability, Theorem~B
    applies.
  \end{proof}

  \noindent
  \parbox[t]{\textwidth}{%
  Kai-Uwe~Bux\\
  Fakul\"at f\"ur Mathematik\\
  Universit\"at Bielefeld\\
  Universit\"atsstra\ss{}e \change[1]{25}\\
  33501 Bielefeld\\
  Germany\\
  website: \texttt{www.kubux.net}
  }
  \bigskip

  \noindent
  \parbox[t]{\textwidth}{%
  Ralf~Gramlich\\
  Fachbereich Mathematik\\
  {\small TU} Darmstadt\\
  Schlo\ss{}gartenstra\ss{}e 7\\
  64289 Darmstadt\\
  Germany\\
  website: \texttt{www.mathematik.tu-darmstadt.de/$\sim$gramlich}\\[1mm]
  School of Mathematics\\
  University of Birmingham\\
  Edgbaston\\
  Birmingham B15 2TT\\
  United Kingdom\\
  website: \texttt{www.mat.bham.ac.uk/staff/gramlichr.shtml}
  }
  \bigskip

  \noindent
  \parbox[t]{\textwidth}{%
  Stefan~Witzel\\
  Fachbereich Mathematik\\
  {\small TU} Darmstadt\\
  Schlo\ss{}gartenstra\ss{}e 7\\
  64289 Darmstadt\\
  Germany\\
  website: \texttt{www.mathematik.tu-darmstadt.de/$\sim$switzel}
  }


\section*{References}
  \begin{references}%
    \newcommand{\nametie}{\,}%
  \begin{article}{Abel91}{Abels:1991}%
    \au{H.\nametie Abels}
    \ti{Finiteness Properties of certain arithmetic groups
        in the function field case}
    \lo{Israel J.\ Math.\ 76 (1991)}{113 -- 128}
  \end{article}
  \begin{thesis}{Abra87}{Abramenko:1987}
    \au{P.\nametie Abramenko}
    \ti{Endlichkeitseigenschaften der Gruppen $\SL_n(\FFF_q[t])$}
    \lo{Thesis Frankfurt am Main (1987)}{}
  \end{thesis}
  \begin{book}{Abra96}{Abramenko:1996}%
    \au{P.\nametie Abramenko}
    \ti{Twin Buildings and Applications to $S$-Arithmetic Groups}
    \lo{Springer LNM~1641 (1996)}{}
  \end{book}
  \begin{book}{AbBr08}{Abramenko.Brown:2008}
    \au{P.\nametie Abramenko, K.S.\nametie Brown}
    \ti{Buildings: Theory and Applications}
    \lo{Springer GTM 248}{}
  \end{book}
  \begin{article}{Behr68}{Behr:1968}%
    \au{H.\nametie Behr}
    \ti{{\selectlanguage{german}Zur starken Approximation in algebraischen Gruppen \"u{}ber globalen
        K\"o{}rpern}}
    \lo{{{\selectlanguage{german}Journal f\"u{}r die reine und angewandte
        Mathematik (Crelles Journal) 229 (1968)}}}{107 -- 116}
  \end{article}
  \begin{article}{Behr98}{Behr:1998}%
    \au{H.\nametie Behr}
    \ti{Arithmetic groups over function fields I\null.
        A complete
        characterization of finitely generated and finitely presented
        arithmetic subgroups of reductive algebraic groups}
    \lo{{\selectlanguage{german}Journal f\"u{}r die reine und
        angewandte Mathematik (Crelles Journal) 495 (1998)}}{79--118}
  \end{article}
  \begin{article}{Behr04}{Behr:2004}
    \au{H.\nametie Behr}
    \ti{Higher finiteness properties of $S$-arithmetic
        groups in the function field case I}
    \lo{In T.W.\,M\"uller (ed.), ``Groups: Geometric and
        Combinatorial Aspects'', Proceedings, Bielefeld 1999;
        LMS LNS 311 (2004)}{27--42}
  \end{article}
  \begin{article}{BoSe76}{Borel.Serre:1976}%
    \au{A.\nametie Borel, J.-P.\nametie Serre}
    \ti{{Cohomologie
    d'immeubles et de groupes S-arithm\'{e}tiques}}
    \lo{Topology 15 (1976)}{211 -- 232}
  \end{article}
  \begin{article}{Bro87}{Brown:1987}%
    \au{K.S.\nametie Brown}
    \ti{Finiteness Properties of Groups}
    \lo{Journal of Pure and Applied Algebra 44 (1987)}{45 -- 75}
  \end{article}
  \begin{book}{Bro89}{Brown:1989}%
    \au{K.S.\nametie Brown}
    \ti{Buildings}
    \lo{Springer Monographs (1989)}{}
  \end{book}
  \begin{article}{BW07}{Bux.Wortman:2007}
    \au{K.-U.\nametie Bux, K.\nametie Wortman}
    \ti{Finiteness properties of arithmetic groups over function fields}
    \lo{Inventiones mathematicae 167 (2007)}{355 -- 378}
  \end{article}
  \begin{article}{BW08}{Bux.Wortman:2008}
    \au{K.-U.\nametie Bux, K.\nametie Wortman}
    \ti{Connectivity Properties of Horospheres in
        Euclidean Buildings and Applications to
        Finiteness Properties of Discrete Groups}
    \lo{Preprint (2008) arXiv:0808.2087}{}
  \end{article}
  \begin{article}{Chev55}{Chevalley:1955}
    \au{C.\nametie Chevalley}
    \ti{Sur certains groupes simples}
    \lo{T\^ohoku Mathematical Journal (1955) 2nd Series 7}{14--66}
  \end{article}
  \begin{article}{v.He03}{Heydebreck:2003}
    \au{A.\nametie von Heydebreck}
    \ti{Homotopy properties of certain complexes associated to sperical
        buildings}
    \lo{Israel Journal of Mathematics 133 (2003)}{369--379}
  \end{article}
  \begin{book}{Kac90}{Kac:1990}
    \au{V.\nametie Kac}
    \ti{Infinite dimensional Lie algebras (third edition)}
    \lo{Cambridge University Press (1990)}{}
  \end{book}
  \begin{book}{Marg91}{Margulis:1991}
    \au{G.A.\nametie Margulis}
    \ti{Discrete Subgroups of Semisimple Lie Groups}
    \lo{Ergebnisse der Mathematik 3.\ Folge, vol.\ 17 (Springer, 1991)}{}
  \end{book}
  \begin{article}{MuVM09}{Maldeghem.Muehlherr:2009}
    \au{B.\nametie M\"uhlherr, H.\nametie Van~Maldeghem}
    \ti{Codistances in buildings}
    \lo{to appear: Innovations in Incidence Geometry}{}
  \end{article}
  \begin{book}{PrSe86}{Pressley.Segal:1986}
    \au{A.\nametie Pressley, G.\nametie Segal}
    \ti{Loop Groups}
    \lo{Clarendon Press, Oxford (1986)}{}
  \end{book}
  \begin{book}{R\'emy02}{Remy:2002}
    \au{{B.\nametie R\'emy}}
    \ti{Groupes de Kac-Moody d\'eploy\'es et presque d\'eploy\'es}
    \lo{Ast\'erisque 277 (2002)}{}
  \end{book}
  \begin{thesis}{Rous77}{Rousseau:1977}
    \au{G.\nametie Rousseau}
    \ti{Immeubles des groupes r\'eductives sur les corps locaux}
    \lo{Thesis, Paris (1977)}{}
  \end{thesis}
  \begin{thesis}{Schu05}{Schulz:2005}
    \au{B.\nametie Schulz}
    \ti{{\selectlanguage{german}Sph"arische Unterkomplexe sph"arischer Geb"aude}}
    \lo{PhD thesis (Frankfurt am Main, 2005)}{}
  \end{thesis}
  \begin{book}{Spri98}{Springer:1998}
    \au{T.A.\nametie Springer}
    \ti{Linear Algebraic Groups (second edition)}
    \lo{Birkh\"auser, Boston (1998)}{}
  \end{book}
  \begin{article}{Stuh80}{Stuhler:1980}%
    \au{U.\nametie Stuhler}
    \ti{Homological properties of certain arithmetic groups in
        the function field case}
    \lo{Inventiones mathematicae 57 (1980)}{263 -- 281}
  \end{article}
  \begin{article}{Tits87}{Tits:1987}%
    \au{J.\nametie Tits}
    \ti{Uniqueness and Presentations of Kac-Moody Groups over Fields}
    \lo{Journal of Algebra~105 (1987)}{542 -- 573}
  \end{article}
  \begin{book}{Weis08}{Weiss:2008}
    \au{R.\nametie Weiss}
    \ti{The structure of affine buildings}
    \lo{Princeton University Press (2008)}{}
  \end{book}
  \begin{article}{????}{????}
    \au{???}
    \ti{???}
    \lo{???}{??? -- ???}
  \end{article}
  \end{references}
\end{document}